\documentclass[leqno,draft]{article}

\def\Hom{\mathop{\rm Hom}}

\newtheorem{theorem}{Theorem}
\newtheorem{lemma}[theorem]{Lemma}
\newtheorem{proposition}[theorem]{Proposition}
\newtheorem{definition}[theorem]{Definition}
\newtheorem{corollary}[theorem]{Corollary}

\newcommand{\begintheorem}{\addtocounter{equation}{1}\begin{theorem}}
\newcommand{\beginlemma}{\addtocounter{equation}{1}\begin{lemma}}
\newcommand{\beginproposition}{\addtocounter{equation}{1}\begin{proposition}}
\newcommand{\begindefinition}{\addtocounter{equation}{1}\begin{definition}}
\newcommand{\begincorollary}{\addtocounter{equation}{1}\begin{corollary}}

\begin{document}

\title{$p$-Adic Heisenberg Cantor sets, 2}

\author{Stephen Semmes \\
        Rice University}

\date{}

\maketitle

\begin{abstract}
In these informal notes, we continue to explore $p$-adic versions of
Heisenberg groups and some of their variants, including the structure
of the corresponding Cantor sets.
\end{abstract}

\tableofcontents

\part{A broad view}
\label{broad view}

\section{Abelian groups}
\label{abelian groups}
\setcounter{equation}{0}

        Let $A$ be an abelian group, in which the group structure is
written additively, and let $A_0 = A \supseteq A_1 \supseteq A_2
\supseteq \cdots$ be a decreasing chain of subgroups of $A$ such that
$\bigcap_{j = 0}^\infty A_j = \{0\}$.  As a basic example, one can
take $A$ to be the group ${\bf Z}$ under addition, $m$ to be a
positive integer greater than or equal to $2$, and $A_j = m^j \, {\bf
Z}$ for each $j \ge 0$, the subgroup of ${\bf Z}$ consisting of
multiples of $m^j$.

        There is a standard way in which to define a topology on $A$
under these conditions, where a set $U \subseteq A$ is an open set if
for each $x \in U$ we have that
\begin{equation}
        x + A_j \subseteq U
\end{equation}
for some $j \ge 0$.  It is easy to see that this does define a
topology on $A$, and that $A$ is Hausdorff with respect to this
topology, because of the condition that $\bigcap_{j = 0}^\infty A_j =
\{0\}$.  One can also check that the group operations of addition and
taking inverses are continuous with respect to this topology, so that
$A$ becomes a topological group.

        If $x \in A$, then let $j(x)$ be the largest nonnegative integer
such that $x \in A_{j(x)}$ when $x \ne 0$, and put $j(0) = + \infty$.  Thus
\begin{equation}
        j(-x) = j(x)
\end{equation}
and
\begin{equation}
        j(x + y) \ge \min(j(x), j(y))
\end{equation}
for every $x, y \in A$.  The second assertion reflects the fact that
$x + y \in A_j$ when $x, y \in A_j$ for some $j \ge 0$.

        Let $r_0 \ge r_1 \ge r_2 \ge \cdots$ be a sequence of positive
real numbers that converges to $0$, and put
\begin{equation}
\label{rho(x) = r_{j(x)}}
        \rho(x) = r_{j(x)}
\end{equation}
when $x \in A$ and $x \ne 0$, and $\rho(0) = 0$.  Observe that
\begin{equation}
\label{rho(-x) = rho(x)}
        \rho(-x) = \rho(x)
\end{equation}
and
\begin{equation}
\label{rho(x + y) le max(rho(x), rho(y))}
        \rho(x + y) \le \max(\rho(x), \rho(y))
\end{equation}
for every $x, y \in A$, by the corresponding properties of $j(x)$.

        If $x, y \in A$, then put
\begin{equation}
\label{d(x, y) = rho(x - y)}
        d(x, y) = \rho(x - y).
\end{equation}
This defines an ultrametric on $A$, which means that is a metric on
$A$ that satisfies
\begin{equation}
\label{d(x, z) le max(d(x, y), d(y, z))}
        d(x, z) \le \max(d(x, y), d(y, z))
\end{equation}
for every $x, y, z \in A$.  Of course, this is a stronger version of
the usual triangle inequality.  By construction, $d(x, y)$ is also
invariant under translations on $A$, in the sense that
\begin{equation}
\label{d(x - z, y - z) = d(x, y)}
        d(x - z, y - z) = d(x, y)
\end{equation}
for every $x, y, z \in A$.  It is easy to see that the topology on $A$
determined by this ultrametic is the same as the one described
earlier.

        Observe that every open or closed ball in $A$ centered at $0$
with positive radius with respect to this ultrametric $d(x, y)$ is
equal to $A_j$ for some $j \ge 0$.  This is a bit nicer when the
$r_j$'s are strictly decreasing, in which case $A_j$ can be expressed
as a ball around $0$ for each $j$.  The situation is a bit better
still when the $A_j$'s are strictly decreasing in $j$, which is to say
that $A_j \ne A_{j + 1}$ for each $j$.

        Conversely, suppose that $d(x, y)$ is a translation-invariant
ultrametric on $A$, and put
\begin{equation}
\label{rho(x) = d(x, 0)}
        \rho(x) = d(x, 0)
\end{equation}
for each $x \in A$.  Thus $\rho(x)$ is a nonnegative real-valued
function on $A$ which is equal to $0$ if and only if $x = 0$.  Because
of translation invariance, (\ref{d(x, y) = rho(x - y)}) holds for
every $x, y \in A$, and one can get (\ref{rho(-x) = rho(x)}) and
(\ref{rho(x + y) le max(rho(x), rho(y))}) from the symmetry of $d(x,
y)$ and the ultrametric version of the triangle inequality for $d(x,
y)$.  Using the ultrametric version of the triangle inequality, we
also get that open and closed balls in $A$ centered at $0$ and with
positive radii with respect to $d(x, y)$ are subgroups of $A$.  A
sequence of these subgroups corresponding to a decreasing sequence of
radii converging to $0$ determines the same topology on $A$ as $d(x,
y)$ does, in the same way as before.

        Note that open balls in any ultrametric space are both open
and closed in the associated topology, and that closed balls are both
open and closed as well.  This is easy to check, using the ultrametric
version of the triangle inequality.  An open subgroup of a topological
group is also closed, because its complement can be expressed as a
union of translates of itself, and hence is an open set too.  If $A$
is equipped with the topology determined by the subgroups $A_j$ as
discussed at the beginning of the section, then each $A_j$ is
automatically an open subgroup of $A$, and thus closed as well.

\section{Non-abelian groups}
\label{non-abelian groups}
\setcounter{equation}{0}

        Let $G$ be a group which is not necessarily abelian, with the
group structure written multiplicatively, and with $e$ as the identity
element.  As in the previous section, let $G_0 = G \supseteq G_1
\supseteq G_2 \supseteq \cdots$ be a decreasing chain of subgroups of
$G$ such that $\bigcap_{j = 0}^\infty = \{e\}$.  Let us say that a set
$U \subseteq G$ is an open set if for each $x \in U$ there is a $j \ge
0$ such that
\begin{equation}
\label{x G_j subseteq U}
        x \, G_j \subseteq U.
\end{equation}
It is easy to see that this defines a topology on $G$, and that this
topology is Hausdorff, because $\bigcap_{j = 0}^\infty G_j = \{e\}$.

        If $x \in G$ and $x \ne e$, then let $j(x)$ be the largest
nonnegative integer such that $x \in G_{j(x)}$, and put $j(e) =
+\infty$.  Observe that
\begin{equation}
\label{j(x^{-1}) = j(x)}
        j(x^{-1}) = j(x)
\end{equation}
and that
\begin{equation}
\label{j(x y) ge min(j(x), j(y))}
        j(x \, y) \ge \min(j(x), j(y))
\end{equation}
for every $x, y \in G$, because the $G_j$'s are subgroups of $G$.  Let
$r_0 \ge r_1 \ge r_2 \ge \cdots$ be a sequence of positive real
numbers that converges to $0$, and put
\begin{equation}
        \rho(x) = r_{j(x)}
\end{equation}
for every $x \in G$ such that $x \ne e$, and $\rho(e) = 0$.  Thus
\begin{equation}
\label{rho(x^{-1}) = rho(x)}
        \rho(x^{-1}) = \rho(x)
\end{equation}
and
\begin{equation}
\label{rho(x y) le max(rho(x), rho(y))}
        \rho(x \, y) \le \max(\rho(x), \rho(y))
\end{equation}
for every $x, y \in G$, by the corresponding properties of $j(x)$.

        If $x, y \in G$, then put
\begin{equation}
\label{d(x, y) = rho(y^{-1} x)}
        d(x, y) = \rho(y^{-1} \, x).
\end{equation}
It is easy to see that this defines an ultrametric on $G$, which
defines the same topology on $G$ as before.  This ultrametric is also
invariant under left translations on $G$, in the sense that
\begin{equation}
\label{d(a x, a y) = d(x, y)}
        d(a \, x, a \, y) = d(x, y)
\end{equation}
for every $a, x, y \in G$.  In particular, $x \mapsto a \, x$ defines
a homeomorphism on $G$ for each $a \in G$.  Of course, one can check
directly that left translations determine homeomorphisms on $G$ with
respect to the topology described earlier.

        Suppose for the moment that the $G_j$'s are normal subgroups
of $G$, so that
\begin{equation}
\label{j(a x a^{-1}) = j(x)}
        j(a \, x \, a^{-1}) = j(x)
\end{equation}
for every $a, x \in G$.  This implies that
\begin{equation}
\label{rho(a x a^{-1}) = rho(x)}
        \rho(a \, x \, a^{-1}) = \rho(x)
\end{equation}
for every $a, x \in G$, and hence that
\begin{equation}
\label{d(x a, y a) = d(x, y)}
        d(x \, a, y \, a) = d(x, y)
\end{equation}
for every $a, x, y \in G$.  In particular, $x \mapsto x \, a$ is a
homeomorphism on $G$ for each $a \in G$ in this case, which could also
be verified more directly.

        Let us say that the chain of subgroups $G_j$ is weakly normal
in $G$ if for each $a \in G$ and $j \ge 0$ there is an $l \ge 0$ such
that
\begin{equation}
\label{G_l subseteq a G_j a^{-1}}
        G_l \subseteq a \, G_j \, a^{-1}.
\end{equation}
This implies that $x \mapsto x \, a$ is a homeomorphism on $G$ for
each $a \in G$, even if the ultrametric $d(x, y)$ may not be invariant
under right translations.  Alternatively, one can first check that $x
\mapsto x^{-1}$ is continuous on $G$ in this case, and then use this
to derive the continuity of right translations from the continuity of
left translations.  One can also check that $G$ is a topological group
under these conditions.

        Actually, the group operations on $G$ are automatically
continuous at the identity element, without this extra hypothesis.
Note that $G_j$ is an open set in $G$ for each $j \ge 0$, by
construction.  Similarly, the left coset $a \, G_j$ of $G_j$ is an
open set in $G$ for each $a \in G$ and $j \ge 0$.  As in the previous
section, $G_j$ is also a closed set in $G$ for each $j$, because its
complement is open.

        Conversely, suppose that $d(x, y)$ is an ultrametric on $G$
that is invariant under left translations, and put
\begin{equation}
\label{rho(x) = d(x, e)}
        \rho(x) = d(x, e)
\end{equation}
for each $x \in G$.  Thus $\rho(x)$ is a nonnegative real-valued
function on $G$ which is equal to $0$ only when $x = 0$.  Because
$d(x, y)$ is invariant under left translations, (\ref{d(x, y) =
rho(y^{-1} x)}) holds for every $x, y \in G$, with this definition of
$\rho$.  Similarly, (\ref{rho(x^{-1}) = rho(x)}) holds for every $x
\in G$, because
\begin{equation}
        \rho(x^{-1}) = d(x^{-1}, e) = d(e, x^{-1}) = d(x, e) = \rho(x).
\end{equation}
This uses the symmetry of $d(\cdot, \cdot)$ in the second step, and
invariance under left translations in the third.  If $x, y \in G$, then
\begin{equation}
        \rho(x \, y) = d(x \, y, e) \le \max(d(x \, y, x), d(x, e))
\end{equation}
by the ultrametric version of the triangle inequality.  This implies that
\begin{equation}
        \rho(x \, y) \le \max(d(y, e), d(x, e)) = \max(\rho(y), \rho(x)),
\end{equation}
using invariance under left translations again.  This shows that
(\ref{rho(x y) le max(rho(x), rho(y))}) also holds for every $x, y \in
G$.  If $d(x, y)$ is invariant under right translations too, then it
is easy to see that (\ref{rho(a x a^{-1}) = rho(x)}) holds for every
$a, x \in G$ as well.

        It follows from (\ref{rho(x^{-1}) = rho(x)}) and (\ref{rho(x
y) le max(rho(x), rho(y))}) that open and closed balls in $G$ centered
at $e$ with positive radii with respect to $d(x, y)$ are subgroups of
$G$.  If $d(x, y)$ is also invariant under right translations on $G$,
then (\ref{rho(a x a^{-1}) = rho(x)}) implies that these balls around
$e$ are normal subgroups of $G$.  Using a decreasing sequence of radii
that converges to $0$, one gets a decreasing sequence of subgroups of
$G$, as before.  The topology on $G$ associated to this sequence of
subgroups is then the same as the topology determined by the given
ultrametric $d(x, y)$.

\section{Topological equivalence}
\label{topological equivalence}
\setcounter{equation}{0}

        Let $G$ be a group, and let $G_0 = G \supseteq G_1 \supseteq
G_2 \supseteq \cdots$ be a decreasing chain of subgroups of $G$ such
that $\bigcap_{j = 1}^\infty G_j = \{e\}$, as in the previous section.
Suppose that
\begin{equation}
        \widetilde{G}_0 = G \supseteq \widetilde{G}_1 \supseteq \widetilde{G}_2
                                                       \supseteq \cdots
\end{equation}
is another decreasing chain of subgroups of $G$ such that $\bigcap_{j
= 1}^\infty \widetilde{G}_j = \{e\}$.  Let us say that these two
chains of subgroups of $G$ are topologically equivalent if for each $j
\ge 1$ there is an $l \ge 1$ such that
\begin{equation}
        \widetilde{G}_l \subseteq G_j,
\end{equation}
and similarly for each $k \ge 1$ there is an $n \ge 1$ such that
\begin{equation}
        G_n \subseteq \widetilde{G}_k.
\end{equation}
This is exactly the condition that ensures that the topologies on $G$
associated to these two chains of subgroups as in the previous section
are the same.

        Suppose that the chain of subgroups of $G$ given by the
$G_j$'s is weakly normal, in the sense described in the preceding
section.  This implies that for each $a \in G$ and $j \ge 1$ there is
a $k \ge 1$ such that
\begin{equation}
\label{a G_k a^{-1} subseteq G_j}
        a \, G_k \, a^{-1} \subseteq G_j,
\end{equation}
by applying the previous condition to $a^{-1}$.  Of course,
\begin{equation}
\label{a G_0 a^{-1} = G supseteq a G_1 a^{-1} supseteq cdots}
        a \, G_0 \, a^{-1} = G \supseteq a \, G_1 \, a^{-1} \supseteq
                                    a \, G_2 \, a^{-1} \supseteq \cdots
\end{equation}
is a decreasing chain of subgroups of $G$ such that $\bigcap_{j =
1}^\infty a \, G_j \, a^{-1} = \{e\}$ for each $a \in G$, because of
the corresponding properties of the $G_j$'s.  Thus weak normality may
be reformulated as saying that (\ref{a G_0 a^{-1} = G supseteq a G_1
a^{-1} supseteq cdots}) is topologically equivalent to the original
chain of subgroups of $G$ given by the $G_j$'s for each $a \in G$.

        Suppose now that $G_j$ has finite index in $G$ for each $j$,
and consider
\begin{equation}
        H_j = \bigcap_{a \in G} a \, G_j \, a^{-1}.
\end{equation}
This is automatically a normal subgroup of $G$, which is also a
subgroup of $G_j$.  Equivalently,
\begin{equation}
        H_j = \bigcap_{a \in A_j} a \, G_j \, a^{-1},
\end{equation}
where $A_j$ is any collection of representatives of the left cosets of
$G_j$ in $G$, so that every element of $G$ can be represented as $a \,
x$ for some $a \in A_j$ and $x \in G_j$.  If $G_j$ has finite index in
$G$, then we can take $A_j$ to have only finitely many elements, and
$H_j$ is the intersection of only finitely many conjugates of $G_j$.
If the chain of $G_j$'s is also weakly normal in $G$, then it follows
that the chain of $H_j$'s is topologically equivalent to the chain of
$G_j$'s in $G$.

\section{Cauchy sequences}
\label{cauchy sequences}
\setcounter{equation}{0}

        Let $G$ be a group, and let $G_0 = G \supseteq G_1 \supseteq
G_2 \supseteq \cdots$ be a decreasing chain of subgroups of $G$ such
that $\bigcap_{j = 0}^\infty G_j = \{e\}$, as in Section
\ref{non-abelian groups}.  Let us say that a sequence $\{x_j\}_{j =
1}^\infty$ of elements of $G$ is a Cauchy sequence if for each
positive integer $n$ there is an $L \ge 1$ such that
\begin{equation}
        x_l^{-1} \, x_j \in G_n
\end{equation}
for every $j, l \ge L$.  This is equivalent to saying that $\{x_j\}_{j
= 1}^\infty$ is a Cauchy sequence with respect to any metric $d(x, y)$
on $G$ that determines the same topology on $G$ as in Section
\ref{non-abelian groups} and is invariant under left translations on
$G$, including the ultrametrics discussed previously.  If $\{x_j\}_{j
= 1}^\infty$ converges to an element $x$ of $G$ with respect to the
topology defined in Section \ref{non-abelian groups}, then it is easy
to see that $\{x_j\}_{j = 1}^\infty$ is a Cauchy sequence in $G$.  If
every Cauchy sequence of elements of $G$ converges to an element of
$G$, then we say that $G$ is complete.

        Suppose that $\widetilde{G}_0 = G \supseteq \widetilde{G}_1
\supseteq \widetilde{G}_2 \supseteq \cdots$ is another decreasing
chain of subgroups of $G$ such that $\bigcap_{j = 1}^\infty
\widetilde{G}_j = \{e\}$.  If the chain of $\widetilde{G}_j$'s is
topologically equivalent to the chain of $G_j$'s in $G$, then it
follows that a sequence $\{x_j\}_{j = 1}^\infty$ of elements of $G$ is
a Cauchy sequence with respect to the chain of $G_j$'s if and only if
it is a Cauchy sequence with respect to the chain of
$\widetilde{G}_j$'s.  In particular, $G$ is the complete with respect
to the chain of $G_j$'s if and only if it is complete with respect to
the chain of $\widetilde{G}_j$'s.

        Let us say that a sequence $\{x_j\}_{j = 1}^\infty$ of
elements of $G$ is ``strongly Cauchy'' if
\begin{equation}
\label{x_l in x_j G_j}
        x_l \in x_j \, G_j
\end{equation}
for every $l \ge j$.  It is easy to see that strongly Cauchy sequences
are Cauchy sequences, and that every Cauchy sequence has a subsequence
that is strongly Cauchy.  As usual, if a Cauchy sequence has a
subsequence that converges, then the whole sequence converges to the
same limit.  It follows that $G$ is complete if every strongly Cauchy
sequence in $G$ converges.

        A pair of Cauchy sequences $\{x_j\}_{j = 1}^\infty$,
$\{y_j\}_{j = 1}^\infty$ in $G$ are said to be ``equivalent'' if
$\{y_j^{-1} \, x_j\}_{j = 1}^\infty$ converges to the identity element
$e$ of $G$.  More explicitly, this means that for each $n \ge 1$ there
is an $L \ge 1$ such that
\begin{equation}
\label{y_j^{-1} x_j in G_n}
        y_j^{-1} \, x_j \in G_n
\end{equation}
for every $j \ge L$.  It is easy to see that this defines an
equivalence relation on the set of Cauchy sequences in $G$, and that a
Cauchy sequence converges if and only if it is equivalent to a
constant sequence.  Note that every subsequence of a Cauchy sequence
$\{x_j\}_{j = 1}^\infty$ is automatically a Cauchy sequence, which is
equivalent to $\{x_j\}_{j = 1}^\infty$.  In particular, every Cauchy
sequence in $G$ is equivalent to a strongly Cauchy sequence, by the
remarks in the previous paragraph.

\section{Cartesian products}
\label{cartesian products}
\setcounter{equation}{0}

        Let $X_1, X_2, \ldots$ be a sequence of nonempty sets, and let
$X = \prod_{j = 1}^\infty X_j$ be their Cartesian product.  Thus $X$
consists of the sequences $x = \{x_j\}_{j = 1}^\infty$ such that $x_j
\in X_j$ for each $j$.  If $x, y \in X$ and $x \ne y$, then let $j(x,
y)$ be the largest nonnegative integer such that $x_j = y_j$ when $j
\le j(x, y)$.  Equivalently, $j(x, y) + 1$ is the smallest positive
integer $j$ such that $x_j \ne y_j$.  In particular, $j(x, y) = 0$
when $x_1 \ne y_1$.  If $x = y$, then put $j(x, y) = + \infty$.
Observe that
\begin{equation}
\label{j(x, y) = j(y, x)}
        j(x, y) = j(y, x)
\end{equation}
and
\begin{equation}
\label{j(x, z) ge min(j(x, y), j(y, z))}
        j(x, z) \ge \min(j(x, y), j(y, z))
\end{equation}
for every $x, y, z \in X$.  More precisely, (\ref{j(x, z) ge min(j(x,
y), j(y, z))}) reduces to the statement that $x_j = z_j$ for every $j
\le l$ when $x_j = y_j$ and $y_j = z_j$ for every $j \le l$.

        Let $r_0 \ge r_1 \ge r_2 \ge \cdots$ be a sequence of positive
real numbers that converges to $0$.  If $x, y \in X$, then put
\begin{equation}
        d(x, y) = r_{j(x, y)}
\end{equation}
when $x \ne y$, and of course $d(x, y) = 0$ when $x = y$.  Thus
\begin{equation}
        d(x, z) \le \max(d(x, y), d(y, z))
\end{equation}
for every $x, y, z \in X$, because of (\ref{j(x, z) ge min(j(x, y), j(y, z))}).
It follows that $d(x, y)$ is an ultrametric on $X$, and it is easy to see
that the topology determined on $X$ by $d(x, y)$ is the same as the
product topology, using the discrete topology on $X_j$ for each $j$.
In particular, $X$ is compact with respect to this topology when $X_j$
has only finitely many elements for each $j$.

        Suppose that $x(1) = \{x_j(1)\}_{j = 1}^\infty, x(2) =
\{x_j(2)\}_{j = 1}^\infty, x(3) = \{x_j(3)\}_{j = 1}^\infty, \ldots$
is a sequence of elements of $X$ which is a Cauchy sequence with
respect to the ultrametric $d(x, y)$.  It is easy to see that this
implies that for each positive integer $j$ there is an $x_j \in X_j$
such that $x_j(l) = x_j$ for all sufficiently large $l$.  If $x =
\{x_j\}_{j = 1}^\infty$, then $x \in X$, and $x(1), x(2), x(3),
\ldots$ converges as a sequence of elements of $X$ to $x$ with respect
to the ultrametric defined in the previous paragraph.  Equivalently,
$x(1), x(2), x(3), \ldots$ converges to $x$ with respect to the
product topology on $X$, corresponding to the discrete topology on
each $X_j$.  This shows that $X$ is complete as a metric space with
respect to the ultrametric defined in the previous paragraph.

\section{An embedding}
\label{embedding}
\setcounter{equation}{0}

        Let $G$ be a group, and let $G_0 = G \supseteq G_1 \supseteq
G_2 \supseteq \cdots$ be a decreasing chain of subgroups such that
$\bigcap_{j = 0}^\infty G_j = \{e\}$, as in Section \ref{non-abelian
groups}.  Even if the $G_j$'s are not normal subgroups of $G$, we can
still consider the quotients $G / G_j$ as sets, where $G / G_j$ is the
set of left cosets $a \, G_j$ of $G_j$ in $G$.  Let $q_j$ be the
canonical quotient mapping from $G$ onto $G_j$, which sends each $a
\in G$ to the corresponding coset $a \, G_j$ in $G / G_j$.  Let $X$ be
the Cartesian product $\prod_{j = 1}^\infty (G / G_j)$ of these
quotient spaces, and let $q$ be the mapping from $G$ into $X$ which
sends $a \in G$ to the sequence $q(a) = \{q_j(a)\}_{j = 1}^\infty$ of
corresponding elements of the quotient spaces $G / G_j$.  Of course,
if $G_j$ is a normal subgroup in $G$ for each $j$, then the quotient
$G / G_j$ is a group in a natural way, and the quotient mapping $q_j$
is a group homomorphism.  In this case, we can also consider $X$ as a
group, where the group operations are defined coordinatewise.  This is
the same as the direct product of the groups $G / G_j$, $j \in {\bf
Z}_+$, and the mapping $q$ is then a homomorphism from $G$ into $X$.

        If $a, b \in G$ and $j \ge 1$, then $q_j(a) = q_j(b)$ if and
only if $a \, G_j = b \, G_j$, which is equivalent to $b^{-1} \, a \in
G_j$.  Let $j_G(a)$ be defined for $a \in G$ as in Section
\ref{non-abelian groups}, and let $j_X(x, y)$ be defined for $x, y \in
X$ as in Section \ref{cartesian products}.  Using the previous remark,
it is easy to see that
\begin{equation}
\label{j_G(b^{-1} a) = j_X(q(a), q(b))}
        j_G(b^{-1} \, a) = j_X(q(a), q(b))
\end{equation}
for every $a, b \in G$.  Let $r_0 \ge r_1 \ge r_2 \ge \cdots$ be a
sequences of positive real numbers that converges to $0$, and let
$d_G$ and $d_X$ be the corresponding ultrametrics on $G$ and $X$
defined in Sections \ref{non-abelian groups} and \ref{cartesian
products}, respectively.  It follows from (\ref{j_G(b^{-1} a) =
j_X(q(a), q(b))}) that
\begin{equation}
        d_G(a, b) = d_X(q(a), q(b))
\end{equation}
for every $a, b \in G$.  In particular, $q$ defines a homeomorphism
from $G$ onto $q(G)$ in $X$, with respect to the topology on $q(G)$
induced by the one on $X$.  Of course, one can also check this
directly, using the original description of the topology on $G$
associated to the $G_j$'s, and the product topology on $X$
corresponding to the discrete topology on $G / G_j$ for each $j$.

\section{Coherent sequences}
\label{coherent sequences}
\setcounter{equation}{0}

        Let $G$ be a group, and let $G_0 = G \supseteq G_1 \supseteq
G_2 \supseteq \cdots$ be a decreasing chain of subgroups of $G$ such
that $\bigcap_{j = 1}^\infty G_j = \{e\}$, as in Section
\ref{non-abelian groups}.  If $j \le l$, then there is a canonical
mapping $\theta_{j, l}$ from $G / G_l$ onto $G / G_j$, which sends a
coset $x \, G_l$ in $G / G_l$ to the corresponding coset $x \, G_j$ in
$G / G_j$.  This is well-defined, because $G_l \subseteq G_j$ when $j
\le l$, by hypothesis.  Note that
\begin{equation}
        \theta_{j, l} \circ \theta_{l, n} = \theta_{j, n}
\end{equation}
when $j \le l \le n$.  If the $G_k$'s are all normal subgroups of $G$,
so that $G / G_k$ is a group for each $k$, then $\theta_{j, l}$ is a
group homomorphism from $G / G_l$ onto $G / G_j$ when $j \le l$.

        Let $a = \{a_j\}_{j = 1}^\infty$ be an element of the
Cartesian product $\prod_{j = 1}^\infty (G / G_j)$, so that $a_j \in G
/ G_j$ for each $j$.  If
\begin{equation}
\label{theta_{j, l}(a_l) = a_j}
        \theta_{j, l}(a_l) = a_j
\end{equation}
for every $j \le l$, then $a = \{a_j\}_{j = 1}^\infty$ is said to be a
coherent sequence.  Of course, it suffices to check that
(\ref{theta_{j, l}(a_l) = a_j}) holds with $l = j + 1$ for each $j$.

        Let $q_j$ be the canonical quotient mapping from $G$ onto $G /
G_j$, as before.  If $x \in G$, then $q(x) = \{q_j(x)\}_{j =
1}^\infty$ is clearly a coherent sequence.  If $\{x_j\}_{j =
1}^\infty$ is a sequence of elements of $G$, then $\{q_j(x_j)\}_{j =
1}^\infty$ is a coherent sequence if and only if $\{x_j\}_{j =
1}^\infty$ is strongly Cauchy, in the sense of Section \ref{cauchy
sequences}.  Note that every element $\{a_j\}_{j = 1}^\infty$ of
$\prod_{j = 1}^\infty (G / G_j)$ can be represented as
$\{q_j(x_j)\}_{j = 1}^\infty$ for some sequence $\{x_j\}_{j =
1}^\infty$ of elements of $G$.  Hence every coherent sequence
$\{a_j\}_{j = 1}^\infty$ can be represented as $\{q_j(x_j)\}_{j =
1}^\infty$ for some strongly Cauchy sequence $\{x_j\}_{j = 1}^\infty$
of elements of $G$.

        If $\{x_j\}_{j = 1}^\infty$ and $\{y_j\}_{j = 1}^\infty$ are
strongly Cauchy sequences of elements of $G$ such that
\begin{equation}
\label{q_j(x_j) = q_j(y_j)}
        q_j(x_j) = q_j(y_j)
\end{equation}
for each $j$, then $\{x_j\}_{j = 1}^\infty$ and $\{y_j\}_{j =
1}^\infty$ are equivalent as Cauchy sequences in $G$, as in Section
\ref{cauchy sequences}.  Conversely, suppose that $\{x_j\}_{j =
1}^\infty$ and $\{y_j\}_{j = 1}^\infty$ are equivalent Cauchy
sequences in $G$, so that $\{y_j^{-1} \, x_j\}_{j = 1}^\infty$
converges to $e$ in $G$.  This implies that for each $j \ge 1$ there
is an $L \ge 1$ such that
\begin{equation}
        y_l^{-1} \, x_l \in G_j
\end{equation}
for every $l \ge L$, as in (\ref{y_j^{-1} x_j in G_n}).  Equivalently,
\begin{equation}
        x_l \in y_l \, G_j \quad\hbox{and}\quad y_l \in x_l \, G_j
\end{equation}
for every $l \ge L$, so that
\begin{equation}
        q_j(x_l) = q_j(y_l)
\end{equation}
for every $l \ge L$.  If $\{x_j\}_{j = 1}^\infty$ and $\{y_j\}_{j =
1}^\infty$ are strongly Cauchy sequences, then it follows that
(\ref{q_j(x_j) = q_j(y_j)}) holds for each $j$, because
\begin{equation}
        q_j(x_l) = q_j(x_j) \quad\hbox{and}\quad q_j(y_l) = q_j(y_j)
\end{equation}
when $l \ge j$.

        If $\{x_j\}_{j = 1}^\infty$ is an ordinary Cauchy sequence in
$G$, then for each $n \ge 1$ there is an $L(n) \ge 1$ such that
\begin{equation}
        x_l^{-1} \, x_j \in G_n
\end{equation}
for every $j, l \ge L(n)$, and hence
\begin{equation}
        x_j \in x_l \, G_n \quad\hbox{and}\quad x_l \in x_j \, G_j
\end{equation}
for every $j, l \ge L(n)$.  Equivalently,
\begin{equation}
        q_n(x_j) = q_n(x_l)
\end{equation}
for every $j, l \ge L(n)$.  Put
\begin{equation}
\label{a_n = q_n(x_{L(n)})}
        a_n = q_n(x_{L(n)})
\end{equation}
for each $n \ge 1$, and observe that $\{a_n\}_{n = 1}^\infty$ is a
coherent sequence.  Let $\{y_j\}_{j = 1}^\infty$ be another Cauchy
sequence in $G$, and let $\{b_n\}_{n = 1}^\infty$ be the corresponding
coherent sequence.  It is easy to see that $a_n = b_n$ for each $n$ if
and only if $\{x_j\}_{j = 1}^\infty$ and $\{y_j\}_{j = 1}^\infty$ are
equivalent as Cauchy sequences in $G$, as in the previous paragraph.

        Suppose that $\prod_{j = 1}^\infty (G / G_j)$ is equipped with
the product topology associated to the discrete topology on $G / G_j$
for each $j$, as before.  It is easy to see that the set of all
coherent sequences is closed as a subset of $\prod_{j = 1}^\infty (G /
G_j)$ with this topology.  Let us check that the image $q(G)$ of $G$
under the usual embedding into $\prod_{j = 1}^\infty (G / G_j)$ is
dense in the set of coherent sequences.  Let $a = \{a_j\}_{j =
1}^\infty$ be a coherent sequence, let $n$ be a positive integer, and
let $x(n)$ be an element of $G$ such that
\begin{equation}
\label{q_n(x(n)) = a_n}
        q_n(x(n)) = a_n.
\end{equation}
Under these conditions, we have that
\begin{equation}
\label{q_j(x(n)) = a_j}
        q_j(x(n)) = a_j
\end{equation}
for $j = 1, \ldots, n$, and hence that $q(x(n))$ converges to $a$ in
$\prod_{j = 1}^\infty (G / G_j)$ as $n \to \infty$, as desired.

        If $G$ is complete, then $q(G)$ is equal to the set of
coherent sequences.  To see this, let $a = \{a_j\}_{j = 1}^\infty$ be
a coherent sequence, and let $\{x_j\}_{j = 1}^\infty$ be a sequence of
elements of $G$ such that $q_j(x_j) = a_j$ for each $j$.  As mentioned
earlier, $\{x_j\}_{j = 1}^\infty$ is a strongly Cauchy sequence in
$G$, which converges to an element $x$ of $G$, because $G$ is
complete.  It follows that $\{x_j\}_{j = 1}^\infty$ is equivalent as a
Cauchy sequence to the constant sequence whose terms are equal to $x$,
so that $q(x) = a$, as desired.  Alternatively, one can use the
completeness of $G$ to show that $q(G)$ is a closed set in $\prod_{j =
1}^\infty (G / G_j)$.  This also uses the fact that $q$ is an
isometric embedding of $G$ into $\prod_{j = 1}^\infty (G / G_j)$ with
respect to suitable metrics, as before.  This implies that $q(G)$ is
the same as the set of coherent sequences, by the remarks in the
preceding paragraph.

        Suppose that $G_j$ is a normal subgroup of $G$ for each $j$,
so that $G / G_j$ is a group for each $j$, and the quotient mapping
$q_j$ is a group homomorphism from $G$ onto $G / G_j$ for each $j$.
As in the previous section, the product $\prod_{j = 1}^\infty (G /
G_j)$ is also a group, where the group operations are defined
coordinatewise, and the usual embedding $q$ of $G$ into $\prod_{j =
1}^\infty (G / G_j)$ is also a group homomorphism.  In this case,
$\prod_{j = 1}^\infty (G / G_j)$ is a topological group with respect
to the product topology associated to the discrete topology on each
factor.  Remember that $q$ is a homeomorphism from $G$ onto $q(G)$,
with respect to the topology on $q(G)$ induced by the product topology
on $\prod_{j = 1}^\infty (G / G_j)$.  The closure $\overline{q(G)}$ of
$q(G)$ in $\prod_{j = 1}^\infty (G / G_j)$ is a subgroup of $\prod_{j
= 1}^\infty (G / G_j)$ under these conditions, and a topological group
with respect to the topology induced by the product topology on
$\prod_{j = 1}^\infty (G / G_j)$.

        Of course, $\overline{q(G)}$ is the same as the set of all
coherent sequences, as before.  Similarly, the closure
$\overline{q(G_l)}$ of $q(G_l)$ in $\prod_{j = 1}^\infty (G / G_j)$
consists of the coherent sequences $a = \{a_j\}_{j = 1}^\infty$ such
that $a_j$ corresponds to the identity element in $G / G_j$ for each
$j \le l$.  This is a normal subgroup of $\overline{q(G)}$ for each
$l$, which is relatively open in $\overline{q(G)}$.  Note that
$\bigcap_{l = 1}^\infty \overline{q(G_l)}$ consists of only the
identity element in $\overline{q(G)}$.  The topology on
$\overline{q(G)}$ determined by the decreasing chain of subgroups
$\overline{q(G_l)}$ as in Section \ref{non-abelian groups} is the same
as the topology induced on $\overline{q(G)}$ by the product topology
on $\prod_{j = 1}^\infty (G / G_j)$.

        Using $q$, we get a natural homomorphism from $G / G_l$ into
$\overline{q(G)} / \overline{q(G_l)}$ for each $l$.  This homomorphism
is injective, because the intersection of $\overline{q(G_l)}$ with
$q(G)$ is equal to $q(G_l)$ for each $l$.  One can check that this
homomorphism is surjective for each $l$ as well, and hence an
isomorphism.  This is basically the same as the density of $q(G)$ in
$\overline{q(G)}$ in this topology.

\section{Completions}
\label{completions}
\setcounter{equation}{0}

        Let $G$ be a group, and let $G_0 = G \supseteq G_1 \supseteq
G_2 \supseteq \cdots$ be a decreasing chain of subgroups of $G$ such
that $\bigcap_{j = 1}^\infty G_j = \{e\}$, as before.  The completion
$\widehat{G}$ of $G$ may be defined as the set of equivalence classes
of Cauchy sequences of elements of $G$, using the notion of
equivalence of Cauchy sequences described in Section \ref{cauchy
sequences}.  Of course, this reduces to $G$ itself when $G$ is
complete, and otherwise there is a natural embedding of $G$ into
$\widehat{G}$, that sends every element $g$ of $G$ to the Cauchy
sequence in which every term is equal to $g$.  There is also a natural
one-to-one correspondence between equivalence classes of Cauchy
sequences in $G$ and coherent sequences, as in the previous section,
so that $\widehat{G}$ can be identified with the set of coherent
sequences.  This is especially nice when $G_j$ is a normal subgroup of
$G$ for each $j$, as discussed later in the previous section.

        Even if the $G_j$'s are not normal subgroups of $G$, we can
still consider $G$ as a metric space with respect to an ultrametric
$d_G(\cdot, \cdot)$ as in Section \ref{non-abelian groups}, and this
completion of $G$ is the same as the completion of $G$ as a metric
space.  If $X$ denotes the product $\prod_{j = 1}^\infty (G / G_j)$,
then we have seen that there is an analogous ultrametric $d_X(\cdot,
\cdot)$ on $X$, such that the usual embedding $q$ of $G$ into $X$ is
an isometric embedding.  We have also seen that $X$ is complete as a
metric space with respect to $d_X(\cdot, \cdot)$, so that the
completion of $G$ as a metric space can be identified with the closure
$\overline{q(G)}$ of $q(G)$ in $X$.  The topology on $X$ determined by
$d_X(\cdot, \cdot)$ is the product topology corresponding to the
discrete topology on each factor, and the closure $\overline{q(G)}$ of
$q(G)$ with respect to this topology is the set of all coherent
sequences.

        Let $a$ be an element of $G$, and consider the corresponding
left translation mapping $g \mapsto a \, g$ on $G$.  It is easy to see
that this mapping sends Cauchy sequences of elements of $G$ to Cauchy
sequences, and that it sends equivalent Cauchy sequences to equivalent
Cauchy sequences.  In particular, this follows from the fact that the
ultrametric $d_G(\cdot, \cdot)$ mentioned in the preceding paragraph
is invariant under left translations.  At any rate, this implies that
$g \mapsto a \, g$ has a natural extension to a mapping on the
completion $\widehat{G}$ of $G$ for each $a \in G$.

        There is also an induced mapping on $G / G_j$ for each $j$,
which sends each left coset $g \, G_j$ in $G / G_j$ to $a \, g \,
G_j$.  This leads to a mapping on $X = \prod_{j = 1}^\infty (G /
G_j)$, using the mappings on $G / G_j$ on each coordinate.  The
restriction of this mapping on $X$ to $q(G)$ corresponds exactly to
the left translation $g \mapsto a \, g$ on $G$.  This mapping on $X$
is an isometry with respect to $d_X(\cdot, \cdot)$, and a
homeomorphism with respect to the product topology associated to the
discrete topology on each factor in particular.  Of course, one can
check directly that this mapping sends coherent sequences to coherent
sequences.

\section{Haar measure}
\label{haar measure}
\setcounter{equation}{0}

        Let $G$ be a group, and let $G_0 = G \supseteq G_1 \supseteq
G_2 \supseteq \cdots$ be a decreasing chain of subgroups of $G$ such
that $\bigcap_{j = 1}^\infty G_j = \{e\}$, as usual.  In this section,
let us suppose also that $G_j$ has finite index in $G$ for each $j$.
This implies that $X = \prod_{j = 1}^\infty (G / G_j)$ is compact with
respect to the product topology associated to the discrete topology on
each factor.  In particular, the closure $\overline{q(G)}$ of the
image $q(G)$ of $G$ under the standard embedding $q$ of $G$ in $X$ is
compact.  If $G$ is complete, then it follows that $G$ is compact with
respect to the topology determined by the $G_j$'s.  Alternatively, the
hypothesis that $G_j$ have finite index in $G$ for each $j$ implies
that $G$ is totally bounded, with respect to any left-invariant metric
that determines the same topology, as in Section \ref{non-abelian
groups}.  This is because $G$ is the union of finitely many left
translates of $G_j$ for each $j$.  This implies that the completion of
$G$ as a metric space is compact, and that $G$ is compact when $G$ is
complete.

        As in the previous section, there is a natural action of $G$
on $G / G_j$ for each $j$, corresponding to left translations on $G$.
This leads to an action of $G$ on $X$ coordinatewise, which restricts
to an action on $\overline{q(G)}$.  There is a natural Borel
probability measure on $\overline{q(G)}$, which is invariant under
this action of $G$ on $\overline{q(G)}$.  Let us begin with the
corresponding invariant integral on $\overline{q(G)}$.

        Let $n$ be a positive integer, and let $\pi_n$ be the obvious
coordinate mapping from $X$ onto $G / G_n$, which sends $a =
\{a_j\}_{j = 1}^\infty$ in $X$ to $\pi_n(a) = a_n$.  Also let $A_n$ be
a subset of $\overline{q(G)}$ such that the restriction of $\pi_n$ to
$A_n$ is a one-to-one mapping onto $G / G_n$.  One may as well take
$A_n$ to be a subset of $q(G)$, so that the elements of $A_n$
correspond to representatives of the left cosets of $G_n$ in $G$.
Note that the number of elements of $A_n$ is equal to the number $|G /
G_n|$ of left cosets of $G_n$ in $G$, which is the same as the index
of $G_n$ in $G$.

        Let $f$ be a continuous real or complex-valued function on
$\overline{q(G)}$, and put
\begin{equation}
\label{I(f, A_n) = frac{1}{|G / G_n|} sum_{a in A_n} f(a)}
        I(f, A_n) = \frac{1}{|G / G_n|} \sum_{a \in A_n} f(a),
\end{equation}
which is the average of $f$ over $A_n$.  If $A_1, A_2, A_3, \ldots$ is
a sequence of subsets of $\overline{q(G)}$ as in the previous
paragraph, then we put
\begin{equation}
\label{I(f) = lim_{n to infty} I(f, A_n)}
        I(f) = \lim_{n \to \infty} I(f, A_n).
\end{equation}
This is analogous to the definition of a Riemann integral as a limit
of Riemann sums.  More precisely, the continuity of $f$ on
$\overline{q(G)}$ and the compactness of $\overline{q(G)}$ implies
that $f$ is uniformly continuous with respect to any of the usual
ultrametrics on $X$.  Using this, one can check that the averages
$I(f, A_n)$ form a Cauchy sequence in ${\bf R}$ or ${\bf C}$, as
appropriate, and that the limit (\ref{I(f) = lim_{n to infty} I(f,
A_n)}) does not depend on the choice of the $A_n$'s.

        One can also check that $I(f)$ is invariant under the action
of $G$ on $\overline{q(G)}$ described earlier.  This does not quite
work for $I(f, A_n)$, but instead the action of $G$ has the effect of
changing $A_n$ to another set with analogous properties for each $n$.
This implies that the limit is invariant under the action of $G$,
because it does not depend on the choice of $A_n$'s, as in the
preceding paragraph.

        Of course, $I(f)$ defines a linear functional on the space of
continuous real or complex-valued functions on $\overline{q(G)}$.
This linear functional is clearly nonnegative, in the sense that $I(f)
\ge 0$ when $f$ is a nonnegative real-valued function on
$\overline{q(G)}$.  The Riesz representation theorem then leads to a
nonnegative Borel measure $\mu$ on $\overline{q(G)}$ such that
\begin{equation}
\label{I(f) = int_G f d mu}
        I(f) = \int_G f \, d\mu
\end{equation}
for every continuous function $f$ on $\overline{q(G)}$.  Note that
$\mu(G) = 1$, because $I(f) = 1$ when $f$ is the constant function
equal to $1$ on $\overline{q(G)}$.  Thus $\mu$ is a probability
measure on $\overline{q(G)}$, which is invariant under the natural
action of $G$ on $\overline{q(G)}$, because of the corresponding
property of $I(f)$.

        If $f_l$ is any real or complex-valued function on $G / G_l$
for some positive integer $l$, then $f_l \circ \pi_l$ defines a
continuous function on $\overline{q(G)}$.  In this case, it is easy to
see that $I(f_l \circ \pi_l, A_l)$ reduces to the average of $f_l$ on
$G / G_l$, and in particular does not depend on the choice of the set
$A_l$ as before.  Moreover, $I(f_l \circ \pi_l, A_n)$ also reduces to
the average of $f_l$ over $G / G_l$ for any $n \ge l$ and any choice
of $A_n$ as before.  This is basically because every left coset of
$G_l$ in $G$ is a union of pairwise-disjoint left cosets of $G_n$ in
$G$.  More precisely, every left coset of $G_l$ in $G$ is the
pairwise-disjoint union of $|G_j / G_n|$ left cosets of $G_n$ in $G$
when $n \ge l$, and $|G / G_n|$ is equal to the product of $|G / G_l|$
and $G_l / G_n|$.

        It follows that $I(f_l \circ \pi_l)$ also reduces to the
average of $f_l$ over $G / G_l$, by passing to the limit as $n \to
\infty$.  Left translations on $G$ simply permute the left cosets of
$G_l$ in $G$, which does not affect the average of $f_l$ over $G /
G_l$.  If $f$ is any continuous real or complex-valued function on
$\overline{q(G)}$, then $f$ can be approximated uniformly on
$\overline{q(G)}$ by functions of the form $f_l \circ \pi_l$, where
$f_l$ is a function on $G / G_l$, because of the uniform continuity of
$f$ mentioned earlier.  Of course, $I(f)$ is a bounded linear
functional on the space of continuous functions on $\overline{q(G)}$
with respect to the supremum norm, because
\begin{equation}
        |I(f)| \le \sup_{a \in \overline{q(G)}} |f(a)|
\end{equation}
for every continuous function $f$ on $\overline{q(G)}$.  Thus $I(f)$
is uniquely determined by its restriction to the dense set of
functions of the form $f_l \circ \pi_l$ on $\overline{q(G)}$.

\section{Heisenberg groups}
\label{heisenberg groups}
\setcounter{equation}{0}

        Let $A$ and $A'$ be abelian groups, and let $B(x, y)$ be an
$A'$-valued function of $x, y \in A$ that is additive in each
variable.  More precisely, this means that
\begin{equation}
\label{B(x + w, y) = B(x, y) + B(w, y)}
        B(x + w, y) = B(x, y) + B(w, y)
\end{equation}
and
\begin{equation}
\label{B(x, y + z) = B(x, y) + B(x, z)}
        B(x, y + z) = B(x, y) + B(x, z)
\end{equation}
for every $w, x, y, z \in A$.  Of course, the group operations on $A$
and $A'$ are expressed additively here.  Note that
\begin{equation}
\label{B(x, 0) + B(x, 0) = B(x, 0 + 0) = B(x, 0)}
        B(x, 0) + B(x, 0) = B(x, 0 + 0) = B(x, 0)
\end{equation}
and hence $B(x, 0) = 0$ for every $x \in A$, and similarly $B(0, y) =
0$ for every $y \in A$.

        Let $G = A \times A'$ as a set, and put
\begin{equation}
\label{(x, s) diamond (y, t) = (x + y, s + t + B(x, y))}
        (x, s) \diamond (y, t) = (x + y, s + t + B(x, y))
\end{equation}
for every $(x, s), (y, t) \in G$, so that $\diamond$ defines a binary
operation on $G$.  If $(w, r)$ is another element of $G$, then
\begin{eqnarray}
\lefteqn{((w, r) \diamond (x, s)) \diamond (y, t)} \\
    & = & (w + x, r + s + B(w, x)) \diamond (y, t) \nonumber \\
    & = & (w + x + y, r + s + t + B(w, x) + B(w + x, y)) \nonumber \\
    & = & (w + x + y, r + s + t + B(w, x) + B(w, y) + B(x, y), \nonumber
\end{eqnarray}
and similarly
\begin{eqnarray}
\lefteqn{(w, r) \diamond ((x, s) \diamond (y, t))} \\
    & = & (w, r) \diamond (x + y, s + t + B(x, y)) \nonumber \\
    & = & (w + x + y, r + s + t + B(w, x + y) + B(x, y)) \nonumber \\
    & = & (w + x + y, r + s + t + B(w, x) + B(w, y) + B(x, y)), \nonumber
\end{eqnarray}
so that $\diamond$ is associative on $G$.  Moreover,
\begin{equation}
\label{(x, s) diamond (0, 0) = (0, 0) diamond (x, s) = (x, s)}
        (x, s) \diamond (0, 0) = (0, 0) \diamond (x, s) = (x, s)
\end{equation}
for every $(x, s) \in G$, since $B(x, 0) = B(0, x) = 0$ for every $x
\in A$, as before.  Thus $(0, 0)$ is the identity element of $G$ with
respect to $\diamond$.

        Observe that
\begin{equation}
 (x, s) \diamond (-x, -s + B(x, x)) = (0, B(x, x) + B(x, -x)) = (0, 0),
\end{equation}
and similarly
\begin{equation}
 (-x, -s + B(x, x)) \diamond (x, s) = (0, B(x, x) + B(-x, x)) = (0, 0)
\end{equation}
for every $(x, s) \in G$.  Thus
\begin{equation}
\label{(x, s)^{-1} = (-x, -s + B(x, x))}
        (x, s)^{-1} = (-x, -s + B(x, x))
\end{equation}
is the inverse element of $G$ corresponding to $(x, s)$ with respect
to $\diamond$.  This is a bit simpler when $B(x, x) = 0$ for every $x
\in A$, in which case
\begin{equation}
\label{(x, s)^{-1} = (-x, -s)}
        (x, s)^{-1} = (-x, -s)
\end{equation}
for every $(x, s) \in G$.  At any rate, this shows that $G$ is a group
with respect to $\diamond$, which is not commutative when $B(x, y)$ is
not symmetric in $x$ and $y$.

\section{Chains of subgroups}
\label{chains of subgroups}
\setcounter{equation}{0}

        Let us continue with the same notations and hypotheses as in
the previous section.  Let $A_0 = A \supseteq A_1 \supseteq A_2
\supseteq \cdots$ be a decreasing chain of subgroups of $A$ such that
$\bigcap_{j = 1}^\infty A_j = \{0\}$, and let $A_0' = A' \supseteq
A_1' \supseteq A_2' \supseteq \cdots$ be a decreasing chain of
subgroups of $A'$ such that $\bigcap_{j = 1}^\infty A_j' = \{0\}$.
Suppose that
\begin{equation}
\label{B(x, y) in A_j' for every x, y in A_j}
        B(x, y) \in A_j' \quad\hbox{for every } x, y \in A_j,
\end{equation}
in which case it is easy to see that $G_j = A_j \times A_j'$ is a
subgroup of $G = A \times A'$ with respect to $\diamond$ for each $j$.
This puts us into the same situation as before, with a decreasing
chain $G_0 = G \supseteq G_1 \supseteq G_2 \supseteq \cdots$ of
subgroups of $G$ whose intersection $\bigcap_{j = 1}^\infty G_j$
consists of only the identity element $(0, 0)$ in $G$.

        Observe that
\begin{eqnarray}
\label{((x, s) diamond (y, t)) diamond (x, s)^{-1} = ...}
\lefteqn{((x, s) \diamond (y, t)) \diamond (x, s)^{-1}} \\
 & = & (x + y, s + t + B(x, y)) \diamond (-x, -s + B(x, x)) \nonumber \\
 & = & (y, t + B(x, y) + B(x + y, -x) + B(x, x)) \nonumber \\
 & = & (y, t + B(x, y) - B(y, x)) \nonumber
\end{eqnarray}
for every $(x, s), (y, t) \in G$.  If
\begin{equation}
\label{B(x, y) - B(y, x) in A_j'}
        B(x, y) - B(y, x) \in A_j'
\end{equation}
for every $x \in A$ and $y \in A_j$, then it follows that $G_j$ is a normal
subgroup of $G$ for each $j$.  In particular, this holds when $B(x, y)$ and
$B(y, x)$ are both elements of $A_j$ for every $x \in A$ and $y \in A_j$.

        Suppose instead that for each $j \ge 1$ there is an $l(j) \ge
j$ such that (\ref{B(x, y) - B(y, x) in A_j'}) holds for every $x \in
A$ and $y \in A_{l(j)}$.  In this case,
\begin{equation}
        \widetilde{G}_j = A_{l(j)} \times A_j'
\end{equation}
is a normal subgroup of $G$ for each $j$, and
\begin{equation}
        G_{l(j)} \subseteq \widetilde{G}_j \subseteq G_j
\end{equation}
for each $j$.  We can extend this to $j = 0$ by putting $l(0) = 0$, so that
\begin{equation}
        \widetilde{G}_0 = G \supseteq \widetilde{G}_1 
                             \supseteq \widetilde{G}_2 \supseteq \cdots
\end{equation}
is a decreasing chain of normal subgroups of $G$ that satisfies
$\bigcap_{j = 1}^\infty \widetilde{G}_j = \{(0, 0)\}$.  Under these
conditions, the chain of $G_j$'s is topologically equivalent to the
chain of $\widetilde{G}_j$'s, as in Section \ref{topological
equivalence}.

        As a variant of this type of condition, suppose that for each
$j \ge 1$ there is an $l(j) \ge j$ such that
\begin{equation}
\label{B(w, z) in A_j'}
        B(w, z) \in A_j'
\end{equation}
for every $w \in A$ and $z \in A_{l(j)}$.  Observe that
\begin{eqnarray}
\label{(y, t)^{-1} diamond (x, s) = ... = (x - y, s - t + B(y, y - x))}
        (y, t)^{-1} \diamond (x, s)
                      & = & (-y, -t + B(y, y)) \diamond (x, s) \\
                      & = & (x - y, s - t + B(y, y) - B(x, y)) \nonumber \\
                      & = & (x - y, s - t + B(y, y - x)) \nonumber
\end{eqnarray}
for every $(x, s), (y, t) \in G$.  If $x - y \in A_{l(j)}$ and $s - t
\in A_j'$, then it follows that (\ref{(y, t)^{-1} diamond (x, s) =
... = (x - y, s - t + B(y, y - x))}) is an element of $A_{l(j)} \times
A_j'$.  Conversely, if (\ref{(y, t)^{-1} diamond (x, s) = ... = (x -
y, s - t + B(y, y - x))}) is an element of $A_{l(j)} \times A_j'$,
then $x - y \in A_{l(j)}$, which implies that
\begin{equation}
        B(y, y - x) \in A_j',
\end{equation}
as in (\ref{B(w, z) in A_j'}).  Thus $s - t \in A_j'$, since $s - t -
B(y, y - x) \in A_j'$, by hypothesis.

        Of course, we can consider $A \times A'$ as an abelian group,
with respect to coordinatewise addition, and we can consider
\begin{equation}
\label{A_0 times A_0' = A times A' supseteq A_1 times A_1' supseteq cdots}
        A_0 \times A_0' = A \times A' \supseteq A_1 \times A_1'
                                    \supseteq A_2 \times A_2' \supseteq \cdots
\end{equation}
as a decreasing chain of subgroups of $A \times A'$ such that
$\bigcap_{j = 1}^\infty A_j \times A_j' = \{(0, 0)\}$, as usual.
The remarks in the previous paragraph show that the topology on $A \times A'$
as an abelian group with respect to this chain of subgroups as in Section 
\ref{abelian groups} is the same as the topology on $G = A \times A'$
with $\diamond$ as the group operation and with respect to the chain of
subgroups $G_j = A_j \times A_j'$, under the condition (\ref{B(w, z) in A_j'}).
More precisely, if $A \times A'$ is equipped with a translation-invariant
metric as in Section \ref{abelian groups}, and if $G$ is equipped with a
metric invariant under left translations with respect to $\diamond$ as
in Section \ref{non-abelian groups}, then the remarks in the previous
paragraph imply that the identity mapping is uniformly continuous as a
mapping from the former to the latter, and from the latter to the former,
under these conditions.

        In particular, a sequence $\{(x_k, s_k)\}_{k = 1}^\infty$ of elements
of $A \times A'$ is a Cauchy sequence in $A \times A'$ as an abelian group
with the chain of subgroups $A_j \times A_j'$ if and only if it is a Cauchy
sequence in $G = A \times A'$ with $\diamond$ as the group operation and with
respect to the chain of subgroups $G_j = A_j \times A_j'$ in this case.
Of course, $\{(x_k, s_k)\}_{k = 1}^\infty$ is a Cauchy sequence in $A \times 
A'$ as an abelian group if and only if $\{x_k\}_{k = 1}^\infty$ and 
$\{s_k\}_{k = 1}^\infty$ are Cauchy sequences in $A$ and $A'$, respectively,
as abelian groups with the chains of subgroups given by the $A_j$'s and the 
$A_j'$'s.  If $A$ and $A'$ are complete as abelian groups with respect to these
chains of subgroups, then it follows that $A \times A'$ is complete as an
abelian group with respect to the chain of subgroups $A_j \times A_j'$, and
that $G = A \times A'$ is complete with  $\diamond$ as the group operation
$\diamond$ with respect to the chain of subgroups $G_j = A_j \times A_j'$.

        It is natural to combine these regularity conditions, and ask
that for each $j \ge 1$ there be an $l(j) \ge 1$ such that
\begin{equation}
\label{B(w, z), B(z, w) in A_j'}
        B(w, z), \, B(z, w) \in A_j'
\end{equation}
for every $w \in A$ and $z \in A_{l(j)}$.  This implies that
(\ref{B(x, y) - B(y, x) in A_j'}) holds for every $x \in A$ and $y \in
A_{l(j)}$, as well as including the previous condition that (\ref{B(w,
z) in A_j'}) hold for every $w \in A$ and $z \in A_{l(j)}$.  
Using this condition, one can also check that if $\{x_k\}_{k =
1}^\infty$ and $\{y_k\}_{k = 1}^\infty$ are Cauchy sequences in $A$ as
an abelian group with respect to the chain of subgroups $A_j$, then
$\{B(x_k, y_k)\}_{k = 1}^\infty$ is a Cauchy sequence in $A'$ as an
abelian group with respect to the chain of subgroups $A_j'$.
Similarly, if $\{w_k\}_{k = 1}^\infty$ and $\{z_k\}_{k = 1}^\infty$ are
Cauchy sequences in $A$ that are equivalent to $\{x_k\}_{k = 1}^\infty$
and $\{y_k\}_{k = 1}^\infty$, respectively, then it is easy to see that
$\{B(w_k, z_k)\}_{k = 1}^\infty$ is equivalent to 
$\{B(x_k, y_k)\}_{k = 1}^\infty$ as a Cauchy sequence in $A'$.  This
permits one to extend $B$ to the completion of $A$ with values in the
completion of $A'$, from which one can then get the completion of $G$.

\section{Haar measure, continued}
\label{haar measure, continued}
\setcounter{equation}{0}

        Let us continue to use the same notation and hypotheses as in
Section \ref{heisenberg groups}.  If $\widetilde{A}$ and
$\widetilde{A}'$ are subgroups of finite index in $A$ and $A'$,
respectively, then $\widetilde{A} \times \widetilde{A}'$ is a subgroup
of finite index in $A \times A'$, where $A \times A'$ is an abelian
group with the group operations defined coordinatewise.  More
precisely, one can identify the quotient $(A \times A') /
(\widetilde{A} \times \widetilde{A}')$ with the product of the
quotients $A / \widetilde{A}$ and $A' / \widetilde{A}'$, so that the
index of $\widetilde{A} \times \widetilde{A}'$ in $A \times A'$ is
equal to the product of the indices of $\widetilde{A}$ in $A$ and
$\widetilde{A}'$ in $A'$.

        Suppose that $B(x, y) \in \widetilde{A}'$ for every $x, y \in
\widetilde{A}$, so that $\widetilde{G} = \widetilde{A} \times
\widetilde{A}'$ is a subgroup of $G = A \times A'$ with respect to
$\diamond$.  Note that $\widetilde{G}_1 = \widetilde{A} \times A'$ is
a normal subgroup of $G$, and that $G / \widetilde{G}_1$ is isomorphic
to $A / \widetilde{A}$.  Similarly, $\widetilde{G}$ is a normal
subgroup of $\widetilde{G}_1$, and $\widetilde{G}_1 / \widetilde{G}$
is isomorphic to $A' / \widetilde{A}'$.  Thus $\widetilde{G}$ has
finite index in $G$, equal to the product of the indices of
$\widetilde{A}$ and $\widetilde{A}'$ in $A$ and $A'$, respectively.

        Let $A_0 = A \supseteq A_1 \supseteq A_2 \supseteq \cdots$ and
$A_0' = A' \supseteq A_1' \supseteq A_2' \supseteq \cdots$ be
decreasing chains of subgroups of $A$ and $A'$, respectively, such
that $\bigcap_{j = 1}^\infty A_j = \{0\}$ and $\bigcap_{j = 1}^\infty
A_j' = \{0\}$, as in the previous section.  Suppose that $A_j$ and
$A_j'$ have finite index in $A$ and $A'$ for each $j$, respectively,
and that $A$ and $A'$ are complete with respect to these chains of
subgroups.  Thus $A$ and $A'$ are compact abelian topological groups
with respect to the topologies determined by these chains of
subgroups, and they have nice translation-invariant Borel probability
measures, as in Section \ref{haar measure}.

        Suppose also that (\ref{B(x, y) in A_j' for every x, y in
A_j}) holds, so that $G_j = A_j \times A_j'$ defines a decreasing
chain of subgroups of $G$ with respect to $\diamond$ that satisfies
$\bigcap_{j = 1}^\infty G_j = \{(0, 0)\}$, as before.  Note that $G_j$
has finite index in $G$ for each $j$, because of the corresponding
properties of the $A_j$'s and $A_j'$'s, by the remarks at the
beginning of the section.  In addition, let us ask that for each $j
\ge 1$ there be an $l(j) \ge 1$ such that (\ref{B(w, z), B(z, w) in
A_j'}) holds for every $w \in A$ and $z \in A_{l(j)}$.  Under these
conditions, $G$ is a topological group with respect to the topology
determined by the $G_j$'s as in Section \ref{non-abelian groups}, and
this topology is the same as the one on $A \times A'$ as an abelian
group corresponding to the subgroups $A_j \times A_j'$.  Of course,
the latter is the same as the product topology on $A \times A'$
corresponding to the topologies on $A$ and $A'$ as abelian groups
determined by the subgroups $A_j$ and $A_j'$, respectively.

        In particular, $G$ is compact with respect to the topology
determined by the $G_j$'s in this case, because $A$ and $A'$ are
compact.  Using the translation-invariant Borel probability measures
on $A$ and $A'$ mentioned earlier, one gets a Borel probability
measure on the product $A \times A'$ that is invariant under
translations on $A \times A'$ as an abelian group, with addition
defined coordinatewise.  A standard argument shows that this measure
is also invariant under translations on $G = A \times A'$ with respect
to $\diamond$, on the left and on the right.  This is more easily seen
in terms of integration of continuous functions on $G$, by integrating
first over $A'$, and then over $A$.

\section{Rings}
\label{rings}
\setcounter{equation}{0}

        Let $R$ be a ring, and let $\mathcal{I}_0 = R \supseteq
\mathcal{I}_1 \supseteq \mathcal{I}_2 \supseteq \cdots$ be a
decreasing chain of ideals in $R$ such that
\begin{equation}
\label{bigcap_{j = 1}^infty mathcal{I}_j = {0}}
        \bigcap_{j = 1}^\infty \mathcal{I}_j = \{0\}.
\end{equation}
Although we are especially interested in commutative rings here, one
can also consider non-commutative rings, in which case one should
specify whether the ideals are one-sided or two-sided.  At any rate,
$R$ is in particular an abelian group with respect to addition, and
the $\mathcal{I}_j$'s form a decreasing chain of subgroups of $R$, as
in Section \ref{abelian groups}.  This leads to a topology on $R$, as
before, and it is easy to see that multiplication on $R$ is continuous
with respect to this topology when $R$ is commutative, or at least
when the $\mathcal{I}_j$'s are two-sided ideals in $R$, so that $R$
becomes a topological ring.  If the $\mathcal{I}_j$'s are all left
ideals, then $x \mapsto a \, x$ is a continuous mapping on $R$ for
every $a \in R$, and similarly for the case where the
$\mathcal{I}_j$'s are all right ideals.  The completion $\widehat{R}$
of $R$ can be defined as before, and multiplication on $R$ can be
extended to $\widehat{R}$ to get a ring when the $\mathcal{I}_j$'s are
two-sided ideals in $R$.  If the $\mathcal{I}_j$'s are left ideals in
$R$, then $x \mapsto a \, x$ can be extended to $\widehat{R}$ for each
$a \in R$, and similarly for the case where the $\mathcal{I}_j$'s are
right ideals.  If $R$ is a commutative ring, then $\widehat{R}$ is a
commutative ring as well.

\section{A basic scenario}
\label{basic scenario}
\setcounter{equation}{0}

        Let $R$ be a commutative ring, and let $N$ be a positive
integer.  Thus $A = R^N$ is a commutative group, where addition is
defined coordinatewise.  We can also think of $A$ as a module over
$R$, with respect to coordinatewise multiplication by elements of $R$.
Let $b = \{b_{p, q}\}_{p, q = 1}^N$ be an $N \times N$ matrix with
entries in $R$, and put
\begin{equation}
\label{B(x, y) = sum_{p = 1}^N sum_{q = 1}^N b_{p, q} x_p y_q}
        B(x, y) = \sum_{p = 1}^N \sum_{q = 1}^N b_{p, q} \, x_p \, y_q
\end{equation}
for each $x = (x_1, \ldots, x_N), y = (y_1, \ldots, y_N) \in R^N$.
This takes values in $A' = R$, and satisfies the additivity conditions
(\ref{B(x + w, y) = B(x, y) + B(w, y)}) and (\ref{B(x, y + z) = B(x,
y) + B(x, z)}) in each variable.  This is also $R$-linear in each
variable, in the sense that
\begin{equation}
\label{B(r x, y) = B(x, r y) = r B(x, y)}
        B(r \, x, y) = B(x, r \, y) = r \, B(x, y)
\end{equation}
for every $x, y \in R^N$ and $r \in R$.  Conversely, if $R$ has a
multiplicative identity element, then any mapping from $R^N \times
R^N$ into $R$ with these properties is of the form (\ref{B(x, y) =
sum_{p = 1}^N sum_{q = 1}^N b_{p, q} x_p y_q}).

        This leads to a group structure $\diamond$ on $G = A \times A'
= R^N \times R$, as in Section \ref{heisenberg groups}.  If $r \in R$,
then put
\begin{equation}
\label{delta_r((x, s)) = (r x, r^2 s)}
        \delta_r((x, s)) = (r \, x, r^2 \, s)
\end{equation}
for each $(x, s) \in G$.  It is easy to see that $\delta_r$ defines a
homomorphism from $G$ into itself for each $r \in R$, and that
\begin{equation}
\label{delta_r circ delta_{widetilde{r}} = delta_{r widetilde{r}}}
        \delta_r \circ \delta_{\widetilde{r}} = \delta_{r \, \widetilde{r}}
\end{equation}
for every $r, r' \in R$.

        Let $\mathcal{I}_0 = R \supseteq \mathcal{I}_1 \supseteq
\mathcal{I}_2 \supseteq \cdots$ be a decreasing chain of ideals in $R$
such that $\bigcap_{j = 1}^\infty \mathcal{I}_j = \{0\}$, as in the
previous section.  We also ask that
\begin{equation}
\label{x y in mathcal{I}_{j + l}}
        x \, y \in \mathcal{I}_{j + l}
\end{equation}
for every $x \in \mathcal{I}_j$ and $y \in \mathcal{I}_l$, and for any
$j, l \ge 0$.  Note that this condition holds automatically when
$\mathcal{I}_j$ is the ``$j$th power'' of $\mathcal{I}_1$, which is
the ideal consisting of all finite sums of products of $j$ elements of
$\mathcal{I}_1$.  If $R$ is the ring ${\bf Z}$ of integers, for
instance, and $\mathcal{I}_1$ is the ideal consisting of integers
divisible by some integer $m \ge 2$, then the $j$th power of
$\mathcal{I}_1$ is the ideal consisting of integers divisible by $m^j$
for each $j \ge 2$.  In this case, the completion $\widehat{R}$ of $R$
is the usual $m$-adic completion of ${\bf Z}$, also known as the
$p$-adic completion of ${\bf Z}$ when $m$ is a prime number $p$.

        Of course, $A_j = \mathcal{I}_j^N$ is a subgroup of $A = R^N$
as an abelian group with respect to coordinatewise addition for each
$j$.  The $A_j$'s form a decreasing chain of subgroups of $A$, for
which the corresponding topology is the same as the product topology
associated to the analogous topology on $R$ for each factor.

        Under these conditions, we have that
\begin{equation}
\label{B(x, y) in mathcal{I}_{j + l}}
        B(x, y) \in \mathcal{I}_{j + l}
\end{equation}
for every $x \in \mathcal{I}_j$ and $y \in \mathcal{I}_l$, and for any
$j, l \ge 0$.  This implies that
\begin{equation}
\label{G_j = A_j times mathcal{I}_{2j} = mathcal{I}_j^N times mathcal{I}_{2j}}
 G_j = A_j \times \mathcal{I}_{2 j} = \mathcal{I}_j^N \times \mathcal{I}_{2 j}
\end{equation}
is a subgroup of $G = R^N \times R$ with respect to $\diamond$ for
each $j \ge 0$.  As usual, the $G_j$'s form a decreasing chain of
subgroups of $G$, whose intersection consists of only the identity
element $(0, 0)$ in $G$.

        Alternatively,
\begin{equation}
\label{H_j = A_j times mathcal{I}_j = mathcal{I}_j^N times mathcal{I}_j}
        H_j = A_j \times \mathcal{I}_j = \mathcal{I}_j^N \times \mathcal{I}_j
\end{equation}
is a normal subgroup of $G$ with respect to $\diamond$ for each $j \ge
0$.  The $H_j$'s also form a decreasing chain of subgroups of $G$,
with $\bigcap_{j = 1}^\infty H_j = \{(0, 0)\}$.  Note that
\begin{equation}
        H_{2 j} \subseteq G_j \subseteq H_j
\end{equation}
for each $j$, so that the chains of $G_j$'s and $H_j$'s are topologically
equivalent in $G$, as in Section \ref{topological equivalence}.

        To be more consistent with the notation used earlier, put
$A_j' = \mathcal{I}_{2 j}$ for each $j$, so that $G_j = A_j \times
A_j'$.  In this case, the regularity conditions discussed in Section
\ref{chains of subgroups} hold with $l(j) = 2 j$.  In particular,
$\widetilde{G}_j = A_{2 j} \times A_j'$ is a normal subgroup of $G$
with respect to $\diamond$ for each $j$.  Of course, this is the same
as $H_{2 j}$ in the present notation.

        The fact that $H_j$ is a normal subgroup of $G$ with respect
to $\diamond$ for each $j$ only uses the hypothesis that $\mathcal{I}_j$
be an ideal in $R$ for each $j$, and not (\ref{x y in mathcal{I}_{j +
l}}).  Similarly, (\ref{x y in mathcal{I}_{j + l}}) is not needed for
the compatibility conditions discussed in Section \ref{chains of
subgroups}, to get in particular that the topology on $G$ determined
by the $H_j$'s is the same as the product topology on $G = R^N \times
R$ associated to the topology on $R$ determined by the $\mathcal{I}_j$'s.
Using (\ref{x y in mathcal{I}_{j + l}}), we get the smaller subgroups $G_j$,
which are equivalent topologically, but correspond to a different type
of geometry.  We also get that the $G_j$'s behave well with respect to
dilations, in the sense that
\begin{equation}
\label{delta_r(G_j) subseteq G_{j + l}}
        \delta_r(G_j) \subseteq G_{j + l}
\end{equation}
for every $r \in \mathcal{I}_l$ and any $l \ge 0$.

\part{Additional structure}
\label{additional structure}

\section{Rings and modules}
\label{rings, modules}
\setcounter{equation}{0}

        Let $R$ be a commutative ring, and let $M$ be a module over
$R$.  Thus $M$ is an abelian group with the group operation expressed
additively, and for each $r \in R$ and $x \in M$ there is an element
$r \, x \in M$ which satisfies the usual properties with respect to
addition and multiplication on $R$ and addition on $M$.  These
properties may be summarized by saying that we have a homomorphism
from $R$ into the ring of endomorphisms on $M$ as an abelian group.
If there is a multiplicative identity element $e$ in $R$, then it is
customary to ask that $e \, x = x$ for every $x \in M$, so that
multiplication by $e$ corresponds to the identity mapping on $M$.

        Of course, $R$ may be considered as a module over itself,
using multiplication on $R$ as a ring to define multiplication by
elements of $R$ on $R$ as a module.  Similarly, ideals in $R$ may be
considered as modules over $R$, and the quotient of $R$ by an ideal is
also a module over $R$ in a natural way.  Any abelian group $A$ may be
considered as a module over the ring ${\bf Z}$ of integers, where $n
\, a$ is the sum of $n$ $a$'s when $n$ is a positive integer, $0 \, a
= 0$, and $(-1) \, a = - a$ for every $a \in A$.  If $R$ is a field,
then a module over $R$ is the same as a vector space over $R$.

        Suppose that $M$ and $N$ are modules over the same commutative
ring $R$.  A mapping $f$ from $M$ into $N$ is said to be a module
homomorphism or $R$-linear if $f$ is a homomorphism from $M$ into $N$
as abelian groups, so that
\begin{equation}
\label{f(x + y) = f(x) + f(y)}
        f(x + y) = f(x) + f(y)
\end{equation}
for every $x, y \in M$, and if 
\begin{equation}
\label{f(r x) = r f(x)}
        f(r \, x) = r \, f(x)
\end{equation}
for every $r \in R$ and $x \in M$.  If $R$ is a field, then this is
the same as a linear mapping between vector spaces.  If $f$ and $g$
are homomorphisms from $M$ into $N$ as $R$-modules and $r \in R$, then
it is easy to see that $f + g$ and $r \, f$ are also homomorphisms
from $M$ into $N$, where $f + g$ and $r \, f$ are defined pointiwse in
the usual way.  The set $\Hom (M, N)$ of homomorphisms from $M$ into
$N$ is thus a module over $R$, which may also be denoted $\Hom_R (M,
N)$, to make the choice of ring $R$ explicit.

        Suppose now that $M_1$, $M_2$, and $N$ are modules over $R$.
A mapping $B$ from $M_1 \times M_2$ into $N$ is said to be
$R$-bilinear if it is $R$-linear in each coordinate separately.  More
precisely, this means that $B(x, y)$ is $R$-linear as a function of
$x$ from $M_1$ into $N$ for each $y \in M_2$, and that $B(x, y)$ is
$R$-linear as a function of $y$ from $M_2$ into $N$ for each $x \in
M_1$.  If $R$ is a field, then this is the same as the usual notion of
bilinear mappings for vector spaces.  We shall be especially
interested here in the case where $M_1 = M_2$.

        As a basic class of examples, let $n$ be a positive integer,
and let $M_1 = M_2$ be the free module $R^n$ of $n$-tuples of elements
of $R$, where addition and multiplication by elements of $R$ are
defined coordinatewise on $R^n$.  Let $b = \{b_{j, l}\}_{j, l = 1}^n$
be an $n \times n$ matrix with entries in $R$, and put
\begin{equation}
\label{B(x, y) = sum_{j = 1}^n sum_{l = 1}^n b_{j, l} x_j y_l}
        B(x, y) = \sum_{j = 1}^n \sum_{l = 1}^n b_{j, l} \, x_j \, y_l
\end{equation}
for every $x, y \in R^n$.  It is easy to see that this is $R$-bilinear
as a mapping from $R^n \times R^n$ into $R$.  Conversely, if $R$ has a
multiplicative identity element, then every $R$-bilinear mapping from
$R^n \times R^n$ into $R$ is of this form.

        Let $M$ be a module over a commutative ring $R$.  By a
submodule of $M$ we mean a subgroup $M'$ of $M$ as an abelian group
such that $r \, x \in M'$ for every $x \in M'$ and $r \in R$.  If we
think of $R$ as a module over itself, then submodules of $R$ are the
same as ideals.  If $R$ is a field, so that a module $M$ over $R$ is
the same as a vector space over $R$, then a submodule of $M$ is the
same as a linear subspace of $M$.  If $R$ is any commutative ring, $M$
and $N$ are modules over $R$, and if $f$ is a module homomorphism from
$M$ into $N$, then the kernel of $f$ is a submodule of $M$, and $f$
maps $M$ onto a submodule of $N$.

        Conversely, let $M$ be a module over a commutative ring $R$,
and let $M'$ be a submodule of $M$.  Because $M'$ is a subgroup of $M$
as an abelian group, we can first form the quotient $M / M'$ as an
abelian group.  It is easy to see that multiplication by elements of
$R$ is well-defined on $M / M'$, so that $M / M'$ is also a module
over $R$, and the natural quotient mapping from $M$ onto $M / M'$ is a
module homomorphism.  Of course, the kernel of this quotient mapping
is equal to $M'$ by construction, so that every submodule is the
kernel of a module homomorphism.

\section{Heisenberg groups}
\label{heisenberg groups, part 2}
\setcounter{equation}{0}

        Let $R$ be a commutative ring, let $M$ and $N$ be $R$-modules,
and let $B$ be an $R$-bilinear mapping from $M \times M$ into $R$.
Consider the binary operation $\diamond$ on $G = M \times N$ defined by
\begin{equation}
\label{(x, s) diamond (y, t) = (x + y, s + t + B(x, y)), part 2}
        (x, s) \diamond (y, t) = (x + y, s + t + B(x, y)).
\end{equation}
One can check that this defines a group structure on $G$, with $(0,
0)$ as the identity element, and
\begin{equation}
\label{(x, s)^{-1} = (-x, -s + B(x, x)), part 2}
        (x, s)^{-1} = (-x, -s + B(x, x))
\end{equation}
as the inverse of $(x, s)$.  If $r \in R$, then
\begin{equation}
\label{delta_r((x, s)) = (r x, r^2 s), part 2}
        \delta_r((x, s)) = (r \, x, r^2 \, s)
\end{equation}
defines a mapping from $G$ into itself, which is a group homomorphism
with respect to $\diamond$.  Observe also that
\begin{equation}
        \delta_{r_1} \circ \delta_{r_2} = \delta_{r_1 r_2}
\end{equation}
for every $r_1, r_2 \in R$.  As in the previous section, if $R$ has a
multiplicative identity element $e$, then it is customary to ask that
multiplication by $e$ act by the identity mapping on $M$ and $N$.  In
this case, $\delta_e$ is also the identity mapping on $G$.

        Suppose that $M'$ and $N'$ are submodules of $M$ and $N$,
respectively, such that
\begin{equation}
\label{B(x, y) in N'}
        B(x, y) \in N'
\end{equation}
for every $x, y \in M'$.  Under these conditions, $M' \times N'$ is a
subgroup of $G$ with respect to $\diamond$, and
\begin{equation}
\label{delta_r(M' times N') subseteq M' times N}
        \delta_r(M' \times N') \subseteq M' \times N'
\end{equation}
for every $r \in R$.  In order to get a normal subgroup of $G$ with
respect to $\diamond$, we should also ask that
\begin{equation}
\label{B(x, y) - B(y, x) in N'}
        B(x, y) - B(y, x) \in N'
\end{equation}
for every $x \in M$ and $y \in M'$.  This is because
\begin{eqnarray}
\label{((x, s) diamond (y, t)) diamond (x, s)^{-1} = ..., part 2}
\lefteqn{((x, s) \diamond (y, t)) \diamond (x, s)^{-1}} \\
  & = & (x + y, s + t + B(x, y)) \diamond (-x, -s + B(x, x)) \nonumber \\
  & = & (y, t + B(x, y) + B(x, x) + B(x + y, - x)) \nonumber \\
  & = & (y, t + B(x, y) - B(y, x)) \nonumber
\end{eqnarray}
for every $(x, s), (y, t) \in G$.

        Suppose instead that
\begin{equation}
\label{B(x, y) and B(y, x) in N'}
        B(x, y) \hbox{ and } B(y, x) \in N',
\end{equation}
for every $x \in M$ and $y \in M'$, which implies (\ref{B(x, y) - B(y,
x) in N'}).  If $\widetilde{M} = M / M'$ and $\widetilde{N} = N / N'$
are the corresponding quotient modules, then it is easy to see that
there is an $R$-bilinear mapping $\widetilde{B}$ from $\widetilde{M}
\times \widetilde{M}$ into $\widetilde{N}$ that corresponds to $B$.
This leads to a group structure $\widetilde{\diamond}$ on
$\widetilde{G} = \widetilde{M} \times \widetilde{N}$, as before.  Let
$q_M$ and $q_N$ be the usual quotient mappings from $M$ and $N$ onto
$\widetilde{M}$ and $\widetilde{N}$, respectively, and let $q$ be the
mapping which sends $(x, s)$ in $G$ to $(q_M(x), q_N(s))$ in
$\widetilde{G}$.  By construction, $q$ is a homomorphism from $G$ onto
$\widetilde{G}$ with respect to $\diamond$ and $\widetilde{\diamond}$,
respectively, and the kernel of $q$ is equal to $M' \times N'$.
Because $\widetilde{M}$ and $\widetilde{N}$ are $R$-modules, we can
define $\widetilde{\delta}_r$ on $\widetilde{G}$ in the same way that
$\delta_r$ was defined on $G$, and with analogous properties.  This
homomorphism $q : G \to \widetilde{G}$ also intertwines $\delta_r$
and $\widetilde{\delta_r}$ for each $r \in R$, in the sense that
\begin{equation}
        q \circ \delta_r = \widetilde{\delta}_r \circ q.
\end{equation}

        Alternatively, suppose that $M_1$ and $N_1$ are $R$-modules,
and that $\phi : M \to M_1$ and $\psi : N \to N_1$ are module
homomorphisms.  Suppose also that $B_1$ is an $R$-bilinear mapping
from $M_1 \times M_1$ into $N_1$, and that
\begin{equation}
\label{B_1(phi(x), phi(y)) = psi(B(x, y))}
        B_1(\phi(x), \phi(y)) = \psi(B(x, y))
\end{equation}
for every $x, y \in M$.  This leads to a group structure $\diamond_1$
on $G_1 = M_1 \times N_1$, as before, and a semigroup of compatible
dilations on $G_1$ associated to elements of $R$ as in
(\ref{delta_r((x, s)) = (r x, r^2 s), part 2}).  It is easy to see
that the mapping from $(x, s)$ in $G$ to $(\phi(x), \psi(s))$ in $G_1$
is a group homomorphism that intertwines the corresponding families of
dilations.  If $M'$ and $N'$ are the kernels of $\phi$ and $\psi$,
then (\ref{B_1(phi(x), phi(y)) = psi(B(x, y))}) implies that
(\ref{B(x, y) and B(y, x) in N'}) holds for every $x \in M$ and $y \in
M'$.

\section{Rings of fractions}
\label{rings of fractions}
\setcounter{equation}{0}

        Let $R$ be a commutative ring with multiplicative identity
element $e$.  Note that $e$ is allowed to be equal to $0$, in which
case $R = \{0\}$.  If $R$ is an integral domain, which is to say that
$R \ne \{0\}$ and there are no nontrivial zero divisors in $R$, then
the corresponding field of fractions can be defined in a standard way.
There are also versions of this that include the possibility of
nontrivial zero divisors, and which are more precise about the
elements of the ring that are allowed as denominators in the
fractions.  We shall review these matters in this section, following
the treatment in Chapter 3 of \cite{a-m}.

        Let $S$ be a multiplicatively closed subset of $R$, in the
sense that $e \in S$ and $s \, t \in S$ for every $s, t \in S$.
Equivalently, $S$ is a sub-semi-group of $R$ as a semigroup with
respect to multiplication, which includes the multiplicative identity
element $e$.  Consider the relation $\equiv$ defined on $R \times S$
by saying that $(a, s) \equiv (b, t)$ when
\begin{equation}
\label{(a t - b s) v = 0}
        (a \, t - b \, s) \, v = 0
\end{equation}
in $R$ for some $v \in S$.  It is easy to see that this relation is
reflexive and symmetric on $R \times S$, and we would like to check
that it is transitive.  If $(a, s)$, $(b, t)$, and $(c, u)$ are
elements of $R \times S$ satisfy $(a, s) \equiv (b, t)$ and $(b, t)
\equiv (c, u)$, then there are $v, w \in S$ such that (\ref{(a t - b
s) v = 0}) and
\begin{equation}
\label{(b u - c t) w = 0}
        (b \, u - c \, t) \, w = 0
\end{equation}
hold.  Multiplying (\ref{(a t - b s) v = 0}) by $u \, w$ and (\ref{(b
u - c t) w = 0}) by $s \, v$, and then adding the resulting equations,
we get that
\begin{equation}
\label{a t v u w - c t w s v = 0}
        a \, t \, v \, u \, w - c \, t \, w \, s \, v = 0.
\end{equation}
This is the same as
\begin{equation}
\label{(a u - c s) t v w = 0}
        (a \, u - c \, s) \, t \, v \, w = 0,
\end{equation}
which implies that $(a, s) \equiv (c, u)$, as desired, since $t \, v
\, w \in S$.

        Thus $\equiv$ defines an equivalence relation on $R \times S$,
and we let $S^{-1} \, R$ denote the collection of corresponding
equivalence classes in $R \times S$.  If $(a, s) \in R \times S$, then
we let $a / s$ denote the equivalence class in $R \times S$ that
contains $(a, s)$.  One can check that
\begin{equation}
\label{(a/s) + (b/t) = (a t + b s) / (s t)}
        (a/s) + (b/t) = (a \, t + b \, s) / (s \, t)
\end{equation}
and
\begin{equation}
\label{(a / s) (b / t) = (a b) / (s t)}
        (a / s) \, (b / t) = (a \, b) / (s \, t)
\end{equation}
are well-defined on $S^{-1} \, R$, and that $S^{-1} \, R$ becomes a
commutative ring with multiplicative identity element given by $e/e$.
By construction,
\begin{equation}
\label{f(a) = a/e}
        f(a) = a/e
\end{equation}
defines a homomorphism from $R$ into $S^{-1} \, R$, which is not
necessarily injective.

        More precisely, it is easy to see that $f(a) = 0$ for some $a
\in R$ if and only if
\begin{equation}
\label{a s = 0}
        a \, s = 0
\end{equation}
for some $s \in S$.  If $R$ is an integral domain, and if $0 \not\in
S$, then it follows that $f$ is injective as a mapping from $R$ into
$S^{-1} \, R$.  In particular, if $R$ is an integral domain, then $S =
R \backslash \{0\}$ is multiplicatively closed, $S^{-1} \, R$ is the
usual field of fractions associated to $R$, and $f$ is the standard
embedding of $R$ into its ring of fractions.  Similarly, for any $R$
and $S$, $S^{-1} \, R = \{0\}$ if and only if $0 \in S$, since one can
apply the previous criterion to $a = e$.

         An element $x$ of $R$ is said to be a unit in $R$ if there is
a $y \in R$ such that
\begin{equation}
\label{x y = e}
        x \, y = e,
\end{equation}
in which case $y$ is unique and may be denoted $x^{-1}$.  One might
also say that a unit $x$ is invertible in $R$, but one should note
that $0$ is a unit when $R = \{0\}$.  As usual, the collection of
units in $R$ is a group with respect to multiplication.

        If $S$ is a multiplicatively closed subset of $R$ and $f : R
\to S^{-1} \, R$ is as before, then $f(s)$ is a unit in $S^{-1} \, R$
for every $s \in S$.  By construction, every element of $S^{-1} \, R$
can be expressed as $f(a) \, f(s)^{-1}$ for some $a \in R$ and $s \in
S$.  If every element of $S$ is a unit in $R$, then $f$ is an
isomorphism from $R$ onto $S^{-1} \, R$, and of course one did not
really need to construct $S^{-1} \, R$.

        As a basic class of examples, let $x$ be an element of $R$,
and let $S$ be the collection of powers $x^n$ of $x$ for each positive
integer $n$, together with $e$, which corresponds to $n = 0$.  If $x$
is a unit in $R$, then every element of $S$ is a unit in $R$, and
$S^{-1} \, R$ is isomorphic to $R$.  If $R$ is the ring ${\bf Z}$ of
integers and $x \ne 0$, then $S^{-1} \, R$ is the ring of rational
numbers with denominator equal to a power of $x$.  If $S$ is any
multiplicatively closed set of nonzero integers, then $S^{-1} \, {\bf
Z}$ is the ring of rational numbers with denominators in $S$.

\section{Modules of fractions}
\label{modules of fractions}
\setcounter{equation}{0}

        Let $R$ be a commutative ring with multiplicative identity
element $e$, let $S$ be a multiplicatively closed subset of $R$, and
let $S^{-1} \, R$ be the corresponding ring of fractions, as in the
previous section.  If $M$ is a module over $R$, then there is an
analogous construction of $S^{-1} \, M$, as a module over $S^{-1} \,
R$.  More precisely, one can first define a relation $\equiv$ on $M
\times S$, by saying that $(x, s) \equiv (x', s')$ for some $x, x' \in
M$ and $s, s' \in S$ when
\begin{equation}
        t \, (s' \, x - s \, x') = 0
\end{equation}
in $M$ for some $t \in S$.  One can check that this is an equivalence
relation on $M \times S$, in essentially the same way as before.

        Let $S^{-1} \, M$ be the corresponding collection of
equivalence classes in $M \times S$, and let $x / s$ denote the
equivalence class that contains $(x, s)$ for every $x \in M$ and $s
\in S$.  One can define addition of elements of $S^{-1} \, M$ in the
same way as before, as well as multiplication of elements of $S^{-1}
\, M$ by elements of $S^{-1} \, R$, so that $S^{-1} \, M$ becomes a
module over $S^{-1} \, R$.  Note that $x / e = 0$ in $S^{-1} \, M$
for some $x \in M$ if and only if
\begin{equation}
\label{s x = 0}
        s \, x = 0
\end{equation}
in $M$ for some $s \in S$.

        If $M = R$ as a module over itself, for instance, then it is
easy to see that $S^{-1} \, M = S^{-1} \, R$ as a module over itself.
Similarly, if $M = R^n$ for some positive integer $n$, where addition
and scalar multiplication are defined coordinatewise, then $S^{-1} \,
M = (S^{-1} \, R)^n$.

        If $M$ is the direct sum of finitely many modules $M_1,
\ldots, M_n$ over $R$, then $S^{-1} \, M$ is equivalent to the direct
sum of $S^{-1} M_1, \ldots, S^{-1} \, M_n$ as modules over $S^{-1} \,
R$.  This also works for the direct sum of infinitely many modules
over $R$, but one should be careful about the distinction between the
direct sum and the direct product for infinitely many modules.  In the
direct sum, all but finitely many coordinates of any element are equal
to $0$, which permits an element of a direct sum of fractions to be
expressed with a common denominator.

        If $M$ and $N$ are modules over $R$ and $\phi : M \to N$ is
$R$-linear, then there is a mapping $S^{-1} \, \phi$ from $S^{-1} \,
M$ into $S^{-1} \, N$ that sends $x / s$ in $S^{-1} \, M$ to $\phi(x)
/ s$ in $S^{-1} \, N$ for every $x \in M$ and $s \in S$.  One can
check that this is well-defined and $S^{-1} \, R$-linear.

        In particular, if $M'$ is a submodule of $M$, then one can
identify $S^{-1} \, M'$ with a submodule of $S^{-1} \, M$.  One can
also verify that $S^{-1} \, (M / M')$ is isomorphic as a module over
$S^{-1} \, R$ to $S^{-1} \, M / S^{-1} \, M'$.

        Now let $M_1$, $M_2$, and $N$ be modules over $R$, and let
$B$ be an $R$-bilinear mapping from $M_1 \times M_2$ into $N$.  
Under these conditions, we get a mapping $S^{-1} \, B$ from $S^{-1} \,
M_1 \times S^{-1} \, M_2$ into $S^{-1} \, N$, which satisfies
\begin{equation}
\label{(S^{-1} B)(x / s, y / t) = B(x, y) / (s t)}
        (S^{-1} \, B)(x / s, y / t) = B(x, y) / (s \, t)
\end{equation}
for every $x \in M_1$, $y \in M_2$, and $s, t \in S$.  Again one can
check that this is well-defined and $S^{-1} \, R$-bilinear.  If $M_1 =
M_2 = R^n$ for some positive integer $n$, and $N = R$, then every
$R$-linear mapping from $M_1 \times M_2$ into $N$ can be expressed as
in (\ref{B(x, y) = sum_{j = 1}^n sum_{l = 1}^n b_{j, l} x_j y_l}).  In
this case, $S^{-1} \, M_1 = S^{-1} \, M_2 = (S^{-1} \, R)^n$, $S^{-1}
\, N = S^{-1} \, R$, and $S^{-1} \, B$ can be given by an analogous
expression, in which the coefficients $b_{j, l} \in R$ are mapped into
$S^{-1} \, R$ in the usual way.

\section{Heisenberg groups of fractions}
\label{heisenberg groups of fractions}
\setcounter{equation}{0}

        Let $R$ be a commutative ring with multiplicative identity
element $e$, let $S$ be a multiplicative subset of $R$, and let
$S^{-1} \, R$ be the corresponding ring of fractions, as in Section
\ref{rings of fractions}.  Also let $M$ and $N$ be modules over $R$,
and let $S^{-1} \, M$ and $S^{-1} \, N$ be the corresponding modules
of fractions over $S^{-1} \, R$, as in the previous section.  If $B$
is an $R$-bilinear mapping from $M \times M$ into $N$, then there is
an associated $S^{-1} \, R$-bilinear mapping $S^{-1} \, B$ from
$S^{-1} \, M \times S^{-1} \, M$ into $S^{-1} \, N$, as in
(\ref{(S^{-1} B)(x / s, y / t) = B(x, y) / (s t)}).  This leads to a
group structure on $S^{-1} \, G = S^{-1} \, M \times S^{-1} \, N$ as
in Section \ref{heisenberg groups, part 2}.  Let $G = M \times N$ with
the group structure associated to $B$ as in Section \ref{heisenberg
groups, part 2}.  There is a natural homomorphism from $G$ into
$S^{-1} \, G$, defined as follows.  There is a natural mapping from
$M$ into $S^{-1} \, M$, which sends $x \in M$ to $x / e$ in $S^{-1} \,
M$.  There is an analogous mapping from $N$ into $S^{-1} \, N$, which
can be combined with the previous one to give a mapping from $G$ into
$S^{-1} \, G$.  It is easy to see that this mapping is a homomorphism,
because of (\ref{(S^{-1} B)(x / s, y / t) = B(x, y) / (s t)}).  This
mapping also intertwines the dilations $\delta_r$ on $G$ with their
counterparts on $S^{-1} \, G$, using the canonical homomorphism from
$R$ into $S^{-1} \, R$ described in Section \ref{rings of fractions}.

\section{Restriction of scalars}
\label{restriction of scalars}
\setcounter{equation}{0}

        Let $R$ and $R_1$ be commutative rings, and let $f$ be a
homomorphism from $R$ into $R_1$.  If $M_1$ is a module over $R_1$,
then we can also think of $M_1$ as a module over $R$, using $f$.  More
precisely, multiplication by $r \in R$ on $M_1$ as a module over $R$
is defined to be the same as multiplication by $f(r)$ on $M_1$ as a
module over $R_1$.  If $R$ and $R_1$ have multiplicative identity
elements $e$ and $e_1$, respectively, then it is customary to require
that $f(e) = e_1$, and that multiplication by $e_1$ correspond to the
identity mapping on $M_1$ as a module over $R_1$, so that multiplication
by $e$ also corresponds to the identity mapping on $M_1$ as a module over $R$.
This process for converting $M_1$ from a module over $R_1$ into a module
over $R$ using $f$ is known as restriction of scalars.

        If $N_1$ is another module over $R_1$ and $\phi$ is a
homomorphism from $M_1$ into $N_1$ as modules over $R_1$, then $\phi$
is also a homomorphism from $M_1$ into $N_1$ as modules over $R$,
using restriction of scalars.  Similarly, $R_1$-bilinear mappings may
be considered as $R$-bilinear mappings using restriction of scalars.

        Let $R$ be a commutative ring with multiplicative identity
element $e$ again, and let $S$ be a multiplicatively closed subset of
$R$.  As in Section \ref{rings of fractions}, this leads to the
corresponding ring of fractions $S^{-1} \, R$, and a natural
homomorphism $f$ from $R$ into $S^{-1} \, R$ associated to the
construction of $S^{-1} \, R$.  If $M$ is a module over $R$, then
$S^{-1} \, M$ may be defined as a module over $S^{-1} \, R$ as in
Section \ref{modules of fractions}.  We can also think of $S^{-1} \,
M$ as a module over $R$, using $f$ and restriction of scalars.  In
particular, the natural mapping
\begin{equation}
        x \mapsto x / e
\end{equation}
from $M$ into $S^{-1} \, M$ may be considered as a homomorphism
between modules over $R$ in this way.

        If $N$ is another module over $R$, then we can also consider
$S^{-1} \, N$ as a module over $R$ using restriction of scalars.  If
$B$ is an $R$-bilinear mapping from $M \times M$ into $N$, then we get
an $S^{-1} \, R$-bilinear mapping $S^{-1} \, B$ from $S^{-1} \, M
\times S^{-1} \, M$ into $S^{-1} \, N$, as in Section \ref{modules of
fractions}.  This leads to a group structure on $S^{-1} \, G = S^{-1}
\, M \times S^{-1} \, N$, as in Section \ref{heisenberg groups, part
2}.  The natural mappings from $M$ and $N$ into $S^{-1} \, M$ and
$S^{-1} \, N$, respectively, combine to give a natural homomorphism
from $G = M \times N$ with the group structure associated to $B$ as in
Section \ref{heisenberg groups, part 2} into $S^{-1} \, G$, as in the
previous section.  These natural mappings from $M$ and $N$ into
$S^{-1} \, M$ and $S^{-1} \, N$ may be considered as homomorphisms
between $R$-modules using restriction of scalars, as in the preceding
paragraph, so that the homomorphism from $G$ into $S^{-1} \, G$
discussed in the previous section can be seen as a homomorphism of the
type described in Section \ref{heisenberg groups, part 2}.

\section{A class of examples}
\label{a class of examples}
\setcounter{equation}{0}

        Let $R$ be a commutative ring with multiplicative identity
element $e$, and let $\mathcal{I}$ be an ideal in $R$.  It is easy to
see that
\begin{equation}
\label{S = e + mathcal{I} = {e + x : x in mathcal{I}}}
        S = e + \mathcal{I} = \{e + x : x \in \mathcal{I}\}
\end{equation}
is a multiplicatively closed subset of $R$.  Note that $0 \not\in S$
when $\mathcal{I} \ne R$.  Thus $S^{-1} \, R$ and $f : R \to S^{-1} \,
R$ can be defined as in Section \ref{rings of fractions}.

        Remember that $f(a) = 0$ for some $a \in R$ if and only if $a
\, s = 0$ for some $s \in S$.  In the present situation, this means that
\begin{equation}
        a \, (e + x) = 0
\end{equation}
for some $x \in \mathcal{I}$.  Equivalently,
\begin{equation}
\label{a = a (-x)}
        a = a \, (-x),
\end{equation}
and hence
\begin{equation}
\label{a = a (-x)^n}
        a = a \, (-x)^n
\end{equation}
for every positive integer $n$.

        If $\mathcal{I}'$, $\mathcal{I}''$ are ideals in $R$, then
their product $\mathcal{I}' \, \mathcal{I}''$ is the ideal in $R$
consisting of all finite sums of products $x \, y$ of elements $x$ of
$\mathcal{I}'$ and $y$ of $\mathcal{I}''$.  Observe that
\begin{equation}
\label{mathcal{I}' mathcal{I}'' = mathcal{I}'' mathcal{I}'}
        \mathcal{I}' \, \mathcal{I}'' = \mathcal{I}'' \, \mathcal{I}',
\end{equation}
because $R$ is commutative, and that
\begin{equation}
\label{mathcal{I}' mathcal{I}'' subseteq mathcal{I}' cap mathcal{I}''}
 \mathcal{I}' \, \mathcal{I}'' \subseteq \mathcal{I}' \cap \mathcal{I}'',
\end{equation}
because $\mathcal{I}'$ and $\mathcal{I}''$ are ideals in $R$.
Similarly, the product of $n$ ideals in $R$ is the ideal consisting of
all finite sums of products of elements from each of the $n$ ideals,
which can also be obtained by taking the products of the ideals one at
a time using the previous definition.  In particular, the $n$th power
of the ideal $\mathcal{I}$ is defined for each positive integer $n$ as
the product of $n$ $\mathcal{I}$'s.  Equivalently, $\mathcal{I}^n$
can be defined recursively by $\mathcal{I}^1 = \mathcal{I}$ and
\begin{equation}
\label{mathcal{I}^{n + 1} = mathcal{I} mathcal{I}^n}
        \mathcal{I}^{n + 1} = \mathcal{I} \, \mathcal{I}^n
\end{equation}
for $n \ge 1$.

        If $a \in R$ satisfies $f(a) = 0$, as before, then there is an
$x \in \mathcal{I}$ such that (\ref{a = a (-x)^n}) holds for each
positive integer $n$, and hence $a \in \mathcal{I}^n$ for each $n$.
In particular, if $\bigcap_{n = 1}^\infty \mathcal{I}^n = \{0\}$, then
it follows that $f : R \to S^{-1} \, R$ is injective.

\section{Chains of submodules}
\label{chains of submodules}
\setcounter{equation}{0}

        Let $R$ be a commutative ring, and let $M$ be a module over
$R$.  Suppose that $M_0 = M \supseteq M_1 \supseteq M_2 \supseteq
\cdots$ is a decreasing chain of submodules of $M$ such that
$\bigcap_{j = 1}^\infty M_j = \{0\}$.  Under these conditions, there
is a standard way to define a topology on $M$, where a subset $U$ of $M$
is an open set if for each $x \in U$ there is a $j \ge 0$ such that
\begin{equation}
\label{x + M_j = {x + y : y in M_j} subseteq U}
        x + M_j = \{x + y : y \in M_j\} \subseteq U.
\end{equation}
It is easy to see that this defines a topology on $M$, and that $x +
M_j$ is an open set in $M$ for each $x \in M$ and $j \ge 0$.  One can
also check that $M$ is Hausdorff with respect to this topology,
because of the hypothesis that $\bigcap_{j = 1}^\infty M_j = \{0\}$.

        By construction, this topology is invariant under translations
on $M$.  More precisely, $M$ becomes a topological abelian group with
respect to addition.  This means that $(x, y) \mapsto x + y$ is a
continuous mapping from $M \times M$ into $M$, where $M \times M$ is
equipped with the product topology corresponding to the topology on
$M$ just described, and that $x \mapsto -x$ is also continuous as a
mapping from $M$ onto itself.  Similarly, for each $r \in R$, one can
check that $x \mapsto r \, x$ is a continuous mapping on $M$.  Note
that the sequence of $M_j$'s forms a local base for this topology on
$M$ at $0$.

        Let $j(x)$ be the largest nonnegative integer such that $x \in
M_{j(x)}$ for each $x \in M$ with $x \ne 0$, and put $j(0) = +\infty$.  Thus
\begin{equation}
\label{j(-x) = j(x)}
        j(-x) = j(x)
\end{equation}
for every $x \in M$, 
\begin{equation}
\label{j(x + y) ge min (j(x), j(y))}
        j(x + y) \ge \min (j(x), j(y))
\end{equation}
for every $x, y \in M$, and
\begin{equation}
\label{j(r x) ge j(x)}
        j(r \, x) \ge j(x)
\end{equation}
for every $x \in M$ and $r \in R$.  Let $a_0 = 1 \ge a_1 \ge a_2 \ge
\cdots$ be a monotone decreasing sequence of positive real numbers
that converges to $0$, and put
\begin{equation}
\label{rho(x) = a_{j(x)}}
        \rho(x) = a_{j(x)}
\end{equation}
for each $x \in M$ with $x \ne 0$, and $\rho(0) = 0$.  Note that
\begin{equation}
\label{rho(-x) = rho(x), part 2}
        \rho(-x) = \rho(x)
\end{equation}
for every $x \in M$,
\begin{equation}
\label{rho(x + y) le max (rho(x), rho(y))}
        \rho(x + y) \le \max (\rho(x), \rho(y))
\end{equation}
for every $x, y \in M$, and
\begin{equation}
\label{rho(r x) le rho(x)}
        \rho(r \, x) \le \rho(x)
\end{equation}
for every $x \in M$ and $r \in R$, by the corresponding properties of $j(x)$.

        It is easy to see that
\begin{equation}
\label{d(x, y) = rho(x - y), part 2}
        d(x, y) = \rho(x - y)
\end{equation}
defines a metric on $M$ that determines the same topology on $M$ as
described earlier.  More precisely, this is an ultrametric on $M$, in
the sense that
\begin{equation}
\label{d(x, z) le max (d(x, y), d(y, z))}
        d(x, z) \le \max (d(x, y), d(y, z))
\end{equation}
for every $x, y, z \in M$.  By construction, this metric is also
invariant under translations, in the sense that
\begin{equation}
\label{d(x - z, y - z) = d(x, y), part 2}
        d(x - z, y - z) = d(x, y)
\end{equation}
for every $x, y, z \in M$.  We also have that
\begin{equation}
        d(r \, x, r \, y) \le d(x, y)
\end{equation}
for every $x, y \in M$ and $r \in R$, because of (\ref{rho(r x) le
rho(x)}).

\section{Completeness}
\label{completeness}
\setcounter{equation}{0}

        Let $R$ be a commutative ring, let $M$ be a module over $R$, and let
\begin{equation}
\label{M_0 = M supseteq M_1 supseteq M_2 supseteq cdots}
        M_0 = M \supseteq M_1 \supseteq M_2 \supseteq \cdots
\end{equation}
be a decreasing chain of submodules of $M$ such that $\bigcap_{j =
1}^\infty M_j = \{0\}$, as in the previous section.  A sequence
$\{x_j\}_{j = 1}^\infty$ of elements of $M$ is said to be a Cauchy
sequence in $M$ if for each $n \ge 1$ there is an $L(n) \ge 1$ such
that
\begin{equation}
\label{x_j - x_l in M_n}
        x_j - x_l \in M_n
\end{equation}
for every $j, l \ge L(n)$.  This is equivalent to saying that
$\{x_j\}_{j = 1}^\infty$ is a Cauchy sequence with respect to any of
the translation-invariant ultrametrics on $M$ discussed in the
previous section, or with respect to any translation-invariant metric
on $M$ that determines the same topology.  As usual, it is easy to see
that every convergent sequence in $M$ is a Cauchy sequence in $M$.
Conversely, if every Cauchy sequence in $M$ converges to an element of
$M$, then $M$ is said to be complete.

        A pair of Cauchy sequences $\{x_j\}_{j = 1}^\infty$ and
$\{y_j\}_{j = 1}^\infty$ of elements of $M$ is said to be equivalent
if $\{x_j - y_j\}_{j = 1}^\infty$ converges to $0$ in $M$.  One can
check that this defines an equivalence relation on the collection of
all Cauchy sequences of elements of $M$.  In particular, every
subsequence of a Cauchy sequence is also Cauchy sequence, and is
equivalent to the original sequence.

        Let us say that a sequence $\{x_j\}_{j = 1}^\infty$ of
elements of $M$ is ``strongly Cauchy'' if
\begin{equation}
\label{x_j - x_l in M_j}
        x_j - x_l \in M_j
\end{equation}
for every $l \ge j \ge 1$.  Every Cauchy sequence in $M$ has a
subsequence which is strongly Cauchy, and hence every Cauchy sequence
in $M$ is equivalent to a strongly Cauchy sequence.  Observe that two
strongly Cauchy sequences $\{x_j\}_{j = 1}^\infty$ and $\{y_j\}_{j =
1}^\infty$ of elements of $M$ are equivalent as Cauchy sequences if
and only if
\begin{equation}
\label{x_j - y_j in M_j}
        x_j - y_j \in M_j
\end{equation}
for each $j$.  This is the same as saying that $x_j$ and $y_j$ have
the same image in $M / M_j$ for each $j$.

        If $j \le l$, then there is a natural homomorphism $\theta_{j,
l}$ from $M / M_l$ onto $M / M_j$ as modules over $R$, because $M_l
\subseteq M_j$.  The kernel of this homomorphism is equal to $M_l /
M_j$, and we have that
\begin{equation}
\label{theta_{j, l} circ theta_{l, n} = theta_{j, n}}
        \theta_{j, l} \circ \theta_{l, n} = \theta_{j, n}
\end{equation}
when $n \ge l \ge j$.  A sequence $\{\xi_j\}_{j = 1}^\infty$ with
$\xi_j \in M / M_j$ for each $j$ is said to be a coherent sequence if
\begin{equation}
        \theta_{j, l}(\xi_l) = \xi_j
\end{equation}
for every $l \ge j$.  Of course, it suffices to check this with $l = j
+ 1$ for each $j$, because of (\ref{theta_{j, l} circ theta_{l, n} =
theta_{j, n}}).

        Note that every sequence $\{\xi_j\}_{j = 1}^\infty$ with
$\xi_j \in M / M_j$ for each $j \ge 1$ can be expressed as
$\{q_j(x_j)\}_{j = 1}^\infty$ for some sequence $\{x_j\}_{j =
1}^\infty$ of elements of $M$, where $q_j$ is the natural quotient
homomorphism from $M$ onto $M / M_j$ as modules over $R$ for each $j$.
It is easy to see that a sequence $\{x_j\}_{j = 1}^\infty$ of elements
of $M$ is strongly Cauchy if and only if $\{q_j(x_j)\}_{j = 1}^\infty$
is a coherent sequence.  Similarly, two strongly Cauchy sequences
$\{x_j\}_{j = 1}^\infty$ and $\{y_j\}_{j = 1}^\infty$ of elements of
$M$ are equivalent as Cauchy sequences in $M$ if and only if
\begin{equation}
        q_j(x_j) = q_j(y_j)
\end{equation}
for every $j$, as in (\ref{x_j - y_j in M_j}).

        If $\{x_j\}_{j = 1}^\infty$ is any Cauchy sequence of elements
of $M$, then $\{q_n(x_j)\}_{j = 1}^\infty$ is eventually constant in
$M / M_n$ for each $n$.  If $\xi_n$ is the limiting value of
$q_n(x_j)$ as $j \to \infty$ for each $n$, then one can check that
$\{\xi_n\}_{n = 1}^\infty$ is a coherent sequence.  If $\{y_j\}_{j =
1}^\infty$ is another Cauchy sequence of elements of $M$, then
$\{y_j\}_{j = 1}^\infty$ is equivalent to $\{x_j\}_{j = 1}^\infty$ if
and only if the corresponding coherent sequences are the same.  Thus
there is a natural one-to-one correspondence between equivalence
classes of Cauchy sequences and coherent sequences.  This
correspondence is a bit simpler when we restrict our attention to
strongly Cauchy sequences in $M$, as in the preceding paragraph.

        It is customary to define the completion $\widehat{M}$ of $M$
to be the set of equivalence classes of Cauchy sequences in $M$.  With
this definition, there is a natural embedding of $M$ into
$\widehat{M}$, which sends each element $x$ of $M$ to the equivalence
class of Cauchy sequences containing the constant sequence $\{x_j\}_{j
= 1}^\infty$ with $x_j = x$ for each $j$.  Alternatively, the
completion of $M$ can be identified with the set of coherent
sequences, as in the previous paragraphs.  The standard embedding of
$M$ into the completion can then be given by sending each $x \in M$ to
\begin{equation}
\label{q(x) = {q_j(x)}_{j = 1}^infty}
        q(x) = \{q_j(x)\}_{j = 1}^\infty,
\end{equation}
which is easily seen to be a coherent sequence.

        Let $X$ be the Cartesian product $\prod_{j = 1}^\infty (M /
M_j)$, consisting of the sequences $\xi = \{\xi_j\}_{j = 1}^\infty$
with $\xi_j \in M / M_j$ for each $j$.  This can also be considered as
a module over $R$, where addition and multiplication by elements of
$R$ are defined termwise.  The mapping $q$ from $M$ into $X$ given by
(\ref{q(x) = {q_j(x)}_{j = 1}^infty}) is a homomorphism of $M$ into
$X$ as modules over $R$, which is injective because of the hypothesis
that $\bigcap_{j = 1}^\infty M_j = \{0\}$.  The set of coherent
sequences forms a submodule of $X$ that contains $q(M)$, so that the
completion $\widehat{M}$ of $M$ may be considered as a module over
$R$.  This is compatible with the usual extension of addition and
scalar multiplication from $M$ to $\widehat{M}$ in terms of Cauchy
sequences.

        We can also consider $X$ as a topological space, using the
product topology associated to the discrete topology on $M / M_j$ for
each $j$.  It is easy to see that $X$ is a topological abelian group
with respect to addition, and that multiplication by $r$ defines a
continuous mapping on $X$ for each $r \in R$.  One can check that $q$
is a homeomorphism from $M$ onto $q(M)$ with respect to the topology
induced on $q(M)$ by the one on $X$, and that the set of coherent
sequences is a closed subset of $X$.  More precisely, the set of
coherent sequences is the same as the closure $\overline{q(M)}$ of
$q(M)$ in $X$, which we identify with the completion $\widehat{M}$ of
$M$.

        If $n$ is a nonnegative integer, then let $X_n$ be the set of
$\xi \in X$ such that
\begin{equation}
\label{xi_j = 0}
        \xi_j = 0
\end{equation}
for each $j \le n$.  Thus $X_0 = X$, $X_{n + 1} \subseteq X_n$ for
each $n$, $X_n$ is a submodule of $X$ for each $n$, and $\bigcap_{n =
1}^\infty X_n = \{0\}$.  This brings us back to the same situation as
for $M$ in the previous section, and it is easy to see that the
topology on $X$ determined by the $X_n$'s as before is the same as the
product topology associated to the discrete topology on $M / M_j$ for
each $j$.  One can also check that $X$ is automatically complete,
basically because Cauchy sequences of elements of $X$ converges
termwise.

        Let $\widehat{M}_n$ be the closure $\overline{q(M_n)}$ of
$q(M_n)$ in $X$ for each $n \ge 0$, which is the same as the
intersection of $\widehat{M} = \overline{q(M)}$ with $X_n$ for each
$n$.  As before, $\widehat{M}_0 = \widehat{M} = \overline{q(M)}$,
$\widehat{M}_n$ is a submodule of $\widehat{M}$ for each $n$,
$\widehat{M}_{n + 1} \subseteq \widehat{M}_n$ for each $n$, and
$\bigcap_{n = 1}^\infty \widehat{M}_n = \{0\}$.  This brings us back
again to the same situation as for $M$ in the previous section, and it
is easy to see that the topology on $\widehat{M}$ determined by the
$\widehat{M}_n$'s is the same as the one induced by the topology on
$X$ already defined.  Note that $\widehat{M}$ is complete, as it
should be.

        We can also consider translation-invariant ultrametrics on $X$
that determine the same topology on $X$, as in the previous section.
If we use the same sequence of $a_j$'s as for $M$, then $q$ is an
isometric embedding of $M$ into $X$, and $\widehat{M} = \overline{q(M)}$
corresponds naturally to the completion of $M$ as a metric space too.

        If $M / M_j$ has only finitely many elements for each $j$,
then $X$ is compact, and it follows that $\widehat{M} =
\overline{q(M)}$ is compact as well.  In this case, it is easy to see
that $M$ was already totally bounded with respect to the
translation-invariant ultrametrics discussed in the previous section.

\section{Chains of ideals}
\label{chains of ideals}
\setcounter{equation}{0}

        Let $R$ be a commutative ring, and suppose that $\mathcal{I}_0
= R \supseteq \mathcal{I}_1 \supseteq \mathcal{I}_2 \supseteq \cdots$
is a decreasing chain of ideals in $R$ such that $\bigcap_{j =
1}^\infty \mathcal{I}_j = \{0\}$.  In particular, we can think of $R$
as a module over itself, and we can think of the $\mathcal{I}_j$'s as
submodules of $R$, so that everything in the previous two sections can
be used here.

        Of course, we have more structure now, and it is easy to see
for instance that $R$ is a topological ring with respect to the
topology defined earlier.  This means that $R$ is a topological
abelian group with respect to addition, as discussed before, and that
multiplication on $R$ is continuous as a mapping from $R \times R$
into $R$, using the product topology on $R \times R$ associated to the
topology already defined on $R$.  If $j(x)$ is defined for $x \in R$
as in Section \ref{chains of submodules}, then we have that
\begin{equation}
\label{j(x y) ge max (j(x), j(y))}
        j(x \, y) \ge \max (j(x), j(y))
\end{equation}
for every $x, y \in R$, in place of (\ref{j(r x) ge j(x)}).  This
corresponds to the fact that
\begin{equation}
\label{x y in mathcal{I}_{j(x)} cap mathcal{I}_{j(y)}}
        x \, y \in \mathcal{I}_{j(x)} \cap \mathcal{I}_{j(y)}
\end{equation}
for every $x, y \in R \backslash \{0\}$.  Similarly, if $\rho(x)$ is
defined for $x \in R$ as in Section \ref{chains of submodules}, using
a monotone decreasing sequence $\{a_j\}_{j = 1}^\infty$ of positive
real numbers as before, then we have that
\begin{equation}
\label{rho(x y) le min (rho(x), rho(y))}
        \rho(x \, y) \le \min (\rho(x), \rho(y))
\end{equation}
for every $x, y \in R$, instead of (\ref{rho(r x) le rho(x)}).

        The completion $\widehat{R}$ of $R$ can be described as in the
preceding section.  The main difference is that multiplication on $R$
can be extended to $\widehat{R}$, so that $\widehat{R}$ is a
commutative ring as well.  More precisely,
\begin{equation}
        X = \prod_{j = 1}^\infty (R / \mathcal{I}_j)
\end{equation}
is a commutative ring in this case, with respect to coordinatewise
addition and multiplication, and using the standard ring structure on
the quotient $R / \mathcal{I}_j$ for each $j$.  If we identify
$\widehat{R}$ with the closure $\overline{q(R)}$ of the image $q(R)$
of the usual embedding $q$ of $R$ into $X$, which is the same as the
set of coherent sequences, then $\widehat{R}$ becomes a subring of
$X$.  This is compatible with the usual extension of multiplication
from $R$ to $\widehat{R}$ in terms of Cauchy sequences.  Note that $X$
is actually a topological ring with respect to the product topology
associated to the discrete topology on $R / \mathcal{I}_j$ for each
$j$, and hence that $\widehat{R} = \overline{q(R)}$ is also a
topological ring with respect to the induced topology.  These
topologies on $X$ and $\widehat{R}$ can be described in terms of
chains of ideals, as in the previous section.

        Suppose that $R$ is a commutative ring with multiplicative
identity element $e$, and that $\mathcal{I}$ is a proper ideal in $R$.
Put $\mathcal{I}_0 = R$, and let $\mathcal{I}_n$ be the $n$th power
$\mathcal{I}^n$ of $\mathcal{I}$ for each positive integer $n$, as in
Section \ref{a class of examples}.  Thus
\begin{equation}
\label{mathcal{I}_{n + 1} subseteq mathcal{I}_n}
        \mathcal{I}_{n + 1} \subseteq \mathcal{I}_n
\end{equation}
for each $n$ by construction, but $\bigcap_{n = 1}^\infty
\mathcal{I}_n = \{0\}$ is an additional hypothesis.  Let $x \in
\mathcal{I}$ be given, and put
\begin{equation}
\label{s_n = sum_{j = 0}^n x^j}
        s_n = \sum_{j = 0}^n x^j
\end{equation}
for each $n$, where $x^j = e$ when $j = 0$.  It is easy to see that
$\{s_n\}_{n = 1}^\infty$ is a strongly Cauchy sequence in $R$ with
respect to this chain of ideals, since $x^n \in \mathcal{I}_n$ for
each $n$.  We also have that
\begin{equation}
        (e - x) \, s_n = e - x^{n + 1}
\end{equation}
for each $n$, by elementary algebra.  It follows that $e - x$ is
invertible in the completion $\widehat{R}$ of $R$, with inverse equal
to the limit of the $s_n$'s in $\widehat{R}$.  Similarly, if $x$ is an
element of the closure $\widehat{I}$ in $\widehat{R}$, which can be
identified with the closure $\overline{q(\mathcal{I})}$ of
$q(\mathcal{I})$ in $X$, then $e - x$ is also invertible in
$\widehat{R}$.

\section{Ideals and submodules}
\label{ideals, submodules}
\setcounter{equation}{0}

        Let $R$ be a commutative ring, and let $M$ be a module over
$R$.  If $\mathcal{I}$ is an ideal in $R$, then the product
$\mathcal{I} \, M$ is defined to be the collection of all finite sums
of products of the form $a \, x$, where $a \in \mathcal{I}$ and $x \in
M$.  It is easy to see that this is a submodule of $M$.

        Suppose that $\mathcal{I}_0 = R \supseteq \mathcal{I}_1
\supseteq \mathcal{I}_2 \supseteq \cdots$ is a decreasing chain of ideals
in $R$ such that $\bigcap_{j = 1}^\infty \mathcal{I}_j = \{0\}$, and that
$M_0 = M \supseteq M_1 \supseteq M_2 \supseteq \cdots$ is a decreasing
chain of submodules of $M$ such that $\bigcap_{j = 1}^\infty M_j = \{0\}$.
As a simple compatibility condition between these two chains, let us
ask that
\begin{equation}
\label{mathcal{I}_j M subseteq M_j}
        \mathcal{I}_j \, M \subseteq M_j
\end{equation}
for each $j$.  This implies that multiplication of elements of $M$ by
elements of $R$ defines a continuous mapping from $R \times M$ into
$M$ with respect to the topologies on $M$ and $R$ discussed in
Sections \ref{chains of submodules} and \ref{chains of ideals}, and
using the corresponding product topology on $R \times M$.  If $j_M(x)$ and
$j_R(r)$ are defined for $x \in M$ and $r \in R$ as before, then we get that
\begin{equation}
\label{j_M(r x) ge max (j_R(r), j_M(x))}
        j_M(r \, x) \ge \max (j_R(r), j_M(x))
\end{equation}
for every $r \in R$ and $x \in M$, in place of (\ref{j(r x) ge j(x)}).
Similarly, if $\rho_M(x)$ and $\rho_R(r)$ are defined for $x \in M$
and $r \in R$ as before, using the same monotone decreasing sequence
$\{a_j\}_{j = 1}^\infty$ of positive real numbers, then we get that
\begin{equation}
\label{rho_M(r x) le min (rho_R(r), rho_M(x))}
        \rho_M(r \, x) \le \min (\rho_R(r), \rho_M(x))
\end{equation}
for every $r \in R$ and $x \in M$, instead of (\ref{rho(r x) le
rho(x)}).

        Note that (\ref{mathcal{I}_j M subseteq M_j}) is equivalent to
saying that the product of any element of $\mathcal{I}_j$ with any
element of $M / M_j$ as a module over $R$ is equal to $0$.  This
implies that $M / M_j$ may also be considered as a module over $R /
\mathcal{I}_j$ for each $j$.  If
\begin{equation}
        X_M = \prod_{j = 1}^\infty (M / M_j) \quad\hbox{and}\quad
           X_R = \prod_{j = 1}^\infty (R / \mathcal{I}_j)
\end{equation}
are as in the previous sections, then it follows that $X_M$ is a
module over $X_R$ with respect to coordinatewise multiplication.
Moreover, multiplication of elements of $X_M$ by elements of $X_R$
defines a continuous mapping from $X_R \times X_M$ into $X_M$ with
respect to the product topologies on $X_M$ and $X_R$ associated to the
discrete topologies on $M / M_j$ and $R / \mathcal{I}_j$ for each $j$,
and using the corresponding product topology on $X_R \times X_M$.

        In this situation, the completion $\widehat{M}$ of $M$ may be
considered as a module over the completion $\widehat{R}$ of $R$.
The completions $\widehat{M}$ and $\widehat{R}$ of $M$ and $R$ may be
identified with the subsets of $X_M$ and $X_R$ consisting of coherent
sequences, respectively, and one can check that the product of a
coherent sequence in $X_M$ by a coherent sequence in $X_R$ is a
coherent sequence in $X_M$.  As usual, this is compatible with extending
scalar multiplication from $M$ and $R$ to $\widehat{M}$ and $\widehat{R}$
using Cauchy sequences under these conditions.  As before, multiplication
of elements of $\widehat{M}$ by elements of $\widehat{R}$ defines a
continuous mapping from $\widehat{R} \times \widehat{M}$ into
$\widehat{M}$ with respect to the appropriate topologies.

\section{Homomorphisms}
\label{homomorphisms}
\setcounter{equation}{0}

        Let $R$ be a commutative ring, let $M$ and $N$ be modules over
$R$, and let $f$ be a homomorphism from $M$ into $N$.  Suppose also
that $M_0 = M \supseteq M_1 \supseteq M_2 \supseteq \cdots$ and $N_0 =
N \supseteq N_1 \supseteq N_2 \supseteq \cdots$ are descreasing chains
of submodules of $M$ and $N$, respectively, such that $\bigcap_{j =
1}^\infty M_j = \{0\}$ and $\bigcap_{j = 1}^\infty N_j = \{0\}$.  As
in Section \ref{chains of submodules}, these chains of submodules
determine topologies on $M$ and $N$.  If $f : M \to N$ is continuous
at $0$, then for each $j \ge 1$ there is an $l(j) \ge 1$ such that
\begin{equation}
\label{f(M_{l(j)}) subseteq N_j}
        f(M_{l(j)}) \subseteq N_j,
\end{equation}
because $N_j$ is a neighborhood of $0$ in $N$.  Conversely, this
condition implies that $f$ is continuous at every point in $M$,
because $f$ is a homomorphism and the topologies on $M$ and $N$ are
invariant under translations.  More precisely, if $x, y \in M$ and $x
- y \in M_{l(j)}$, then (\ref{f(M_{l(j)}) subseteq N_j}) implies that
\begin{equation}
\label{f(x) - f(y) = f(x - y) in N_j}
        f(x) - f(y) = f(x - y) \in N_j,
\end{equation}
because $f : M \to N$ is a homomorphism.  This is basically a uniform
continuity condition for $f$.  If $M$ and $N$ are equipped with
translation-invariant metrics that determine the same topologies, for
instance, then $f$ is uniformly continuous with respect to these
metrics on $M$ and $N$.  In particular, this implies that $f$ maps
Cauchy sequence in $M$ to Cauchy sequences in $N$, and that $f$ can be
extended to a uniformly continuous mapping $\widehat{f}$ from the
completion $\widehat{M}$ of $M$ into the completion $\widehat{N}$ of
$N$.  It is easy to see that $\widehat{f}$ is also a homomorphism from
$\widehat{M}$ into $\widehat{N}$ as modules over $R$.

        This is all a bit simpler if (\ref{f(M_{l(j)}) subseteq N_j})
holds with $l(j) = j$, so that
\begin{equation}
\label{f(M_j) subseteq N_j}
        f(M_j) \subseteq N_j
\end{equation}
for each $j$.  In this case, $f$ induces a homomorphism $f_j$ from $M
/ M_j$ into $N / N_j$ as modules over $R$ for each $j$.  As usual,
\begin{equation}
        X_M = \prod_{j = 1}^\infty (M / M_j) \quad\hbox{and}\quad
                X_N = \prod_{j = 1}^\infty (N / N_j)
\end{equation}
may be considered as modules over $R$, where addition and
multiplication by elements of $R$ are defined coordinatewise.  Thus
$f$ determines a homomorphism $f_X$ from $X_M$ into $X_N$ as modules
over $R$, using the induced homomorphisms $f_j$ on each factor.  Note
that $f_X$ is also continuous with respect to the product topologies
on $X_M$ and $X_N$ associated to the discrete topologies on $M / M_j$
and $N / N_j$ for each $j$.  If we identify the completions $\widehat{M}$
and $\widehat{N}$ with the subsets of $X_M$ and $X_N$ consisting of
coherent sequences, respectively, then the extension $\widehat{f}$ of
$f$ to a mapping from $\widehat{M}$ into $\widehat{N}$ corresponds to
the restriction of $f_X$ to the set of coherent sequences in $X_M$.
In particular, it is easy to see that $f_X$ sends coherent sequences
in $X_M$ to coherent sequences in $X_M$.

        Suppose in addition that $\mathcal{I}_0 = R \supseteq 
\mathcal{I}_1 \supseteq \mathcal{I}_2 \supseteq \cdots$ is a decreasing
chain of ideals in $R$ such that $\bigcap_{j = 1}^\infty \mathcal{I}_j = 
\{0\}$, and that $\mathcal{I}_j \, M \subseteq M_j$ and $\mathcal{I}_j \, N 
\subseteq N_j$ for each $j$, as in the previous section.  This implies
that $M / M_j$ and $N / N_j$ are also modules over $R / \mathcal{I}_j$
for each $j$, and one can check that $f_j$ is a homomorphism from $M /
M_j$ into $N / N_j$ as modules over $R / \mathcal{I}_j$ for each $j$
too.  As in the previous section, $X_M$ and $X_N$ may be considered as
modules over $X_R = \prod_{j = 1}^\infty (R / \mathcal{I}_j)$, and now
we get that $f_X$ is a homomorphism from $X_M$ into $X_N$ as modules
over $X_R$ as well.  Under these conditions, we get that the extension
$\widehat{f}$ of $f$ to a mapping from the completion $\widehat{M}$ of
$M$ to the completion $\widehat{N}$ of $N$ is a homomorphism from
$\widehat{M}$ into $\widehat{N}$ as modules over the completion
$\widehat{R}$ of $R$.  The same conclusion could also be obtained
using the weaker continuity condition (\ref{f(M_{l(j)}) subseteq N_j})
and Cauchy sequences, and one could weaken the compatibility condition
between the ideals $\mathcal{I}_j$ and the submodules $M_j$ and $N_j$
analogously, but some of the other steps would not work in the same
way.

        As another variant, $M$ and $N$ may be commutative rings
themselves, with decreasing chains of ideals $M_j$ and $N_j$ as
before.  If $f$ is now a continuous homomorphism from $M$ into $N$,
then $f$ extends to a continuous homomorphism $\widehat{f}$ from the
completion $\widehat{M}$ of $M$ to the completion $\widehat{N}$ of $N$
as commutative rings.  This uses the fact that products of Cauchy
sequences in $M$ are still Cauchy sequences in $M$, which are mapped
by $f$ to the corresponding products of Cauchy sequences in $N$.  The
description of $\widehat{f}$ is more elegant when $f$ satisfies the
stronger continuity condition (\ref{f(M_j) subseteq N_j}), so that
$f_j : M / M_j \to N / N_j$ is a ring homomorphism for each $j$.  This
implies that $f_X$ is ring homomorphism from $X_M$ into $X_N$, and
$\widehat{f}$ may be identified with the restriction of $f_X$ to the
set of coherent sequences in $X_M$.

\section{Bilinear mappings}
\label{bilinear mappings}
\setcounter{equation}{0}

        Let $R$ be a commutative ring, let $M$ and $N$ be modules over
$R$, and let $B$ be an $R$-bilinear mapping from $M \times M$ into
$N$.  One could also consider bilinear mappings defined on products of
different modules over $R$, but this is will not be needed here.
Suppose that
\begin{equation}
 M_0 = M \supseteq M_1 \supseteq M_2 \supseteq \cdots \quad\hbox{and}\quad 
  N_0 = N \supseteq N_1 \supseteq N_2 \supseteq \cdots
\end{equation}
are decreasing chains of submodules of $M$ and $N$, respectively, such that 
\begin{equation}
        \bigcap_{j = 1}^\infty M_j = \{0\} \quad\hbox{and}\quad
         \bigcap_{j = 1}^\infty N_j = \{0\},
\end{equation}
as usual.  As a simple compatibility condition, let us ask that
\begin{equation}
\label{B(x, y) in N_j}
        B(x, y) \in N_j
\end{equation}
for every $x \in M$ and $y \in M_j$, and for every $x \in M_j$ and $y
\in M$, for each $j \ge 1$.  This is the same as saying that $B(x, y)$
satisfies the continuity condition (\ref{f(M_j) subseteq N_j}) as a
function of $x$ for each $y \in M$, and also as a function of $y$ for
each $x \in M$.

        It is easy to see that $B$ induces an $R$-bilinear mapping
\begin{equation}
\label{B_j : (M / M_j) times (M / M_j) to N / N_j}
        B_j : (M / M_j) \times (M / M_j) \to N / N_j
\end{equation}
for each $j$, and hence an $R$-bilinear mapping
\begin{equation}
\label{B_X : X_M times X_M to X_N}
        B_X : X_M \times X_M \to X_N,
\end{equation}
where $X_M = \prod_{j = 1}^\infty (M / M_j)$ and $X_N = \prod_{j =
1}^\infty (N / N_j)$, as before.  One can also check that $B_X$ sends
pairs of coherent sequences in $X_M$ to coherent sequences in $X_N$.
If we identify the completions $\widehat{M}$ and $\widehat{N}$ with
the subsets of $X_M$ and $X_N$ consisting of coherent sequences,
respectively, then the restriction of $B_X$ to coherent sequences in
$X_M$ defines a natural extension of $B$ to an $R$-bilinear mapping
$\widehat{B}$ from $\widehat{M} \times \widehat{M}$ into $\widehat{N}$.
This is equivalent to extending $B$ to an $R$-bilinear mapping
$\widehat{B}$ from $\widehat{M} \times \widehat{M}$ into $\widehat{N}$
using Cauchy sequences, which also works under less stringent continuity
conditions.  More precisely, one can use continuity conditions like these
on $B$ to show that $B$ sends pairs of Cauchy sequences in $M$ to Cauchy
sequences in $N$, and that equivalent pairs of Cauchy sequences in $M$
are sent to equivalent Cauchy sequences in $N$, and so on.

        Suppose now that $\mathcal{I}_0 = R \supseteq \mathcal{I}_1
\supseteq \mathcal{I}_2 \supseteq \cdots$ is a decreasing chain of
ideals in $R$ such that $\bigcap_{j = 1}^\infty \mathcal{I}_j =
\{0\}$, and that $\mathcal{I}_j \, M \subseteq M_j$ and $\mathcal{I}_j
\, N \subseteq N_j$ for each $j$, as in the previous sections.  Thus
$M / M_j$ and $N / N_j$ are also modules over $R / \mathcal{I}_j$ for
each $j$, and one can check that $B_j$ is $R/\mathcal{I}_j$-bilinear
for each $j$ too.  As before, $X_M$ and $X_N$ may be considered as
modules over $X_R = \prod_{j = 1}^\infty (R / \mathcal{I}_j)$, and it
is easy to see that $B_X$ is $X_R$-bilinear.  As in Section
\ref{ideals, submodules}, the completions $\widehat{M}$ and
$\widehat{N}$ of $M$ and $N$ may be considered as modules over the
completion $\widehat{R}$ of $R$ in this case, and one can check that
$\widehat{B}$ is $\widehat{R}$-bilinear.  This can be verified more
directly in terms of Cauchy sequences as well, using the fact that the
product of a Cauchy sequence in $M$ or $N$ with a Cauchy sequence in
$R$ is still a Cauchy sequence in $M$ or $N$, respectively, and so on.

        Let $\diamond$ be the group structure on $G = M \times N$
associated to $B$ as in Section \ref{heisenberg groups, part 2}.  The
hypothesis (\ref{B(x, y) in N_j}) implies that $H_j = M_j \times N_j$
is a normal subgroup of $G$ for each $j$, as before.  Similarly, there
is a group structure $\diamond_j$ on
\begin{equation}
\label{G_j = (M / M_j) times (N / N_j)}
        G_j = (M / M_j) \times (N / N_j)
\end{equation}
associated to $B_j$ for each $j$.  Let
\begin{equation}
        \Phi_j : G \to G_j
\end{equation}
be the obvious quotient mapping, obtained by combining the canonical
quotient mappings from $M$ onto $M / M_j$ and $N$ onto $N / N_j$.  As
in Section \ref{heisenberg groups, part 2}, $\Phi_j$ is a homomorphism
from $G$ onto $G_j$ for each $j$, with kernel equal to $H_j$.

\section{Stronger conditions}
\label{stronger conditions}
\setcounter{equation}{0}

        Let us continue with the same hypotheses and notations as in
the previous section.  As a stronger version of (\ref{B(x, y) in
N_j}), let us ask that
\begin{equation}
\label{B(x, y) in N_{j + l}}
        B(x, y) \in N_{j + l}
\end{equation}
for every $x \in M_j$, $y \in M_l$, and $j, l \ge 0$.  This implies that
\begin{equation}
\label{E_j = M_j times N_{2 j}}
        E_j = M_j \times N_{2 j}
\end{equation}
is a subgroup of $G = M \times N$ with respect to $\diamond$ for each
$j \ge 1$.  Note that
\begin{equation}
\label{H_{2 j} subseteq E_j subseteq H_j}
        H_{2 j} \subseteq E_j \subseteq H_j
\end{equation}
for each $j$, but that $E_j$ is not necessarily a normal subgroup of
$G$.

        Suppose that $\mathcal{I}_0 = R \supseteq \mathcal{I}_1
\supseteq \mathcal{I}_2 \supseteq \cdots$ is a decreasing chain of
ideals in $R$ such that $\bigcap_{j = 1}^\infty \mathcal{I}_j = \{0\}$
and
\begin{equation}
\label{mathcal{I}_j mathcal{I}_l subseteq mathcal{I}_{j + l}}
        \mathcal{I}_j \, \mathcal{I}_l \subseteq \mathcal{I}_{j + l}
\end{equation}
for every $j, l \ge 0$.  Note that (\ref{mathcal{I}_j mathcal{I}_l
subseteq mathcal{I}_{j + l}}) holds automatically when $\mathcal{I}_j$
is equal to the $j$th power $\mathcal{I}^j$ of some ideal
$\mathcal{I}$ in $R$ for each $j \ge 1$.  As a stronger version of the
usual compatibility conditions
\begin{equation}
\label{mathcal{I}_j M subseteq M and mathcal{I}_j N subseteq N}
        \mathcal{I}_j \, M \subseteq M \quad\hbox{and}\quad
          \mathcal{I}_j \, N \subseteq N,
\end{equation}
let us ask that
\begin{equation}
\label{I_j M_l subseteq M_{j + l} and I_j N_l subseteq N_{j + l}}
        \mathcal{I}_j \, M_l \subseteq M_{j + l} \quad\hbox{and}\quad
         \mathcal{I}_j \, N_l \subseteq N_{j + l}
\end{equation}
for each $j, l \ge 0$.  This would follow from (\ref{mathcal{I}_j
mathcal{I}_l subseteq mathcal{I}_{j + l}}) if $M_l = \mathcal{I}_l \,
M$ and $N_l = \mathcal{I}_l \, N$ for each $l$, as would (\ref{B(x, y)
in N_{j + l}}).  Under these conditions, we get that
\begin{equation}
\label{delta_r(E_l) subseteq E_{j + l}}
        \delta_r(E_l) \subseteq E_{j + l}
\end{equation}
for every $r \in \mathcal{I}_j$ and $j, l \ge 0$, where $\delta_r$
is as in (\ref{delta_r((x, s)) = (r x, r^2 s), part 2}).

        Let us look at how these various conditions extend to the
corresponding completions.  If $\widehat{M}$ is the completion of $M$,
and $\widehat{M}_j$ is the closure of $M_j$ in $\widehat{M}$, then
$\widehat{M}_j$ is a submodule of $\widehat{M}$ with $\widehat{M}_{j +
1} \subseteq \widehat{M}_j$ for each $j$, and $\bigcap_{j = 1}^\infty
\widehat{M}_j = \{0\}$.  As usual, it is easier to deal with
completions in terms of coherent sequences, although one can also work
directly with Cauchy sequences.  Similar remarks apply to $N$ and the
$N_j$'s, and the extension $\widehat{B}$ of $B$ to an $R$-bilinear
mapping from $\widehat{M} \times \widehat{M}$ into $\widehat{N}$ was
discussed in the previous section.  In the context of the previous
section, $\widehat{B}$ would satisfy the analogue of (\ref{B(x, y) in
N_j}) for $\widehat{M}_j$ and $\widehat{N}_j$, and here $\widehat{B}$
satisfies the analogue of (\ref{B(x, y) in N_{j + l}}) for $\widehat{M}_j$
and $\widehat{N}_j$.  In the same way, if $\widehat{\mathcal{I}}_j$ is the
closure of $\widehat{\mathcal{I}}$ in the completion $\widehat{R}$ of $R$,
then $\widehat{\mathcal{I}}_j$ is an ideal in $\widehat{R}$ such that
$\widehat{\mathcal{I}}_{j + 1} \subseteq \widehat{\mathcal{I}}_j$ for
each $j$, and $\bigcap_{j = 1}^\infty \widehat{\mathcal{I}}_j = \{0\}$.
As before, $\widehat{M}$ and $\widehat{N}$ may be considered as modules
over $\widehat{R}$, and $\widehat{B}$ is $\widehat{R}$-bilinear.  It is easy 
to see that (\ref{mathcal{I}_j mathcal{I}_l subseteq mathcal{I}_{j + l}})
and (\ref{I_j M_l subseteq M_{j + l} and I_j N_l subseteq N_{j + l}}) imply
their counterparts for $\widehat{\mathcal{I}}_j$, $\widehat{M}_j$, and 
$\widehat{N}_j$ as well.  If $\widehat{G} = \widehat{M} \times \widehat{N}$, 
then we can extend the group structure $\diamond$ on $G$ to the group structure
$\widehat{\diamond}$ on $\widehat{G}$ associated to $\widehat{B}$ as in Section
\ref{heisenberg groups, part 2}, and $\widehat{H}_j = \widehat{M}_j \times
\widehat{N}_j$ is a normal subgroup of $\widehat{G}$ for each $j$.
Under the conditions of this section, $\widehat{E}_j = \widehat{M}_j \times
\widehat{N}_{2 j}$ is subgroup of $\widehat{G}$ that satisfies the
counterparts of (\ref{H_{2 j} subseteq E_j subseteq H_j}) and
(\ref{delta_r(E_l) subseteq E_{j + l}}) with $\widehat{H}_j$ and
$\widehat{\mathcal{I}}_j$ instead of $H_j$ and $\mathcal{I}_j$.

        As usual, one can get basic examples by taking $R$ to be the
ring ${\bf Z}$ of integers, and $\mathcal{I}$ to be the ideal $a \,
{\bf Z}$ of integer multiples of some integer $a \ge 2$.  If
$\mathcal{I}_j$ is the $j$th power $\mathcal{I}^j$ of $\mathcal{I}$
when $j \ge 1$, then $\mathcal{I}_j$ is the same as the ideal $a^j \,
{\bf Z}$ of integer multiples of $a^j$.  In this case, the completion
$\widehat{R}$ of $R$ is the same as the $a$-adic completion of ${\bf Z}$.

\end{document}